\magnification=\magstep1
\input amstex
\documentstyle{amsppt}
\TagsOnRight
\hsize=5.25in                                                  
\vsize=7.8in

\nopagenumbers

\magnification=\magstep1

\topmatter

\title Novikov - Shubin signatures, I\endtitle
\rightheadtext{M. Farber}
\leftheadtext{Novikov - Shubin signatures, I}
\author  Michael Farber \endauthor
\address
School of Mathematical Sciences,
Tel-Aviv University,
Ramat-Aviv 69978, Israel
\endaddress
\email farber\@math.tau.ac.il
\endemail
\thanks{The research was  supported by a grant from the Israel Academy of Sciencies
and Humanities and by the Herman Minkowski Center for Geometry}
\endthanks
\abstract{Torsion objects of von Neumann categories describe the phenomenon "spectrum near zero"
discovered by S. Novikov and M. Shubin.
In this paper we classify Hermitian forms on torsion objects of a finite von Neumann category.
We prove that any such Hermitian form can be represented as the discriminant form of 
a degenerate Hermitian form on a projective module. We also find the relation between the Hermitian forms on projective modules which holds if and only if their discriminant forms are congruent.
A notion of superfinite von Neumann category is introduced. 
It is proven that the classification of torsion Hermitian forms in a
superfinite category can be completely reduced to the isomorphism types 
of their positive and the negative parts.}
\endabstract
\endtopmatter

\define\C{{\bold C}}
\define\R{{\bold R}}

\define\T{{\Cal T}}  
\define\X{{\Cal X}}
\define\CClZ{{\Cal Z}}   
\define\Y{{\Cal Y}}

\define\A{{\Cal A}}  
  
\redefine\B{{\Cal B}}  
\define\e{{\frak e}}

\define\sign{\operatorname{sign}} 
 
\define\im{\operatorname{im}} 
   

\define\id{\operatorname{id}}    
\define\Id{\operatorname{Id}}   
     
\define\ob{\operatorname{Ob}}    
\define\cl{\operatorname{cl}} 
\redefine\H{\Cal H}

\define\E{\Cal E}
  
\def\<{\langle}
\def\>{\rangle}

\define\pd#1#2{\dfrac{\partial#1}{\partial#2}}
\documentstyle{amsppt}   
\define\CC{\Cal C}
\define\eca{\Cal E({\Cal C})}     
\define\tca{\T(\CC)}      
\define\cc{\Cal C}

\redefine\L{\Cal L}
\define\lam{\lambda}
\define\CClK{\Cal K}

\define\tdim {{\frak {tordim}}}


\define\Hilb{\frak {Hilb}}
\define\lm{\operatorname{\Lambda-mod}}

\nopagenumbers

\heading{\bf \S 0. Introduction}\endheading

S. Novikov and M. Shubin \cite{NS1}, \cite{NS2} discovered a new way of producing topological
invariants of manifolds by studying the spectrum near zero of the Laplacian acting on $L^2$ forms
on the universal cover. It was proven by M. Gromov and M. Shubin \cite{GS1}, \cite{GS2}
that the Novikov - Shubin invariants depend only on the homotopy type of the manifold.
W. L\"uck and J. Lott \cite{LL} computed the Novikov - Shubin invariants for some 3-dimensional
manifolds.

Two different homological "explanations" the Novikov - Shubin invariants were 
suggested independently by W. L\"uck \cite{Lu1}, \cite{Lu2},\cite{Lu3} and by the author
\cite{Fa1}, \cite{Fa2}. 

The main idea of the approach developed in \cite{Fa1 - Fa2}
was to view the "spectrum near zero" as {\it a torsion part of the $L^2$ cohomology,
understood in an extended sense}. The main principles of \cite{Fa1-Fa2}
are the following: {\it 

$\bullet$ the category of Hilbert representations of a given descrete group can be canonically embedded into an abelian category $\E$;

$\bullet$ the extended abelian category $\E$ contains the chain complex 
$\ell^2(\pi)\otimes_{\pi} C_\ast(\tilde M)$ of $L^2$ forms on the universal covering $\tilde M$
for any compact polyhedron $M$;

$\bullet$ the homology of this chain complex
(called {\it the extended $L^2$ cohomology}) is a well defined object of $\E$. 

$\bullet$ the extended cohomology naturally splits as a sum of its projective and 
torsion parts; 

$\bullet$ the theory of von Neumann dimension provides a tool to measure the size of 
the projective part of the extended cohomology;

$\bullet$ the Novikov-Shubin invariants depend only of the torsion part of the extended cohomology.}

There are numerous advantages of viewing the Novikov - Shubin invariants as byproducts of some
homology theory. Firstly, it immediately implies their homotopy invariance 
(the result of Gromov and Shubin). Secondly,
there exist other invariants of torsion part of the extended $L^2$-homology, 
(for example, the ones constructed in \cite{Fa1}, \cite{Fa2} and others) 
which are independent of the Novikov - Shubin invariants. Thirdly, using some of these
invariants allows strengthening (obtained in \cite{Fa1}, \cite{Fa2}) of the Morse type inequalities
of Novikov and Shubin \cite{NS1}, \cite{NS2}. 
At last, as we showed in \cite{Fa1}, \cite{Fa2} 
the well developed tools of homological algebra 
(such as derived functors and spectral sequences) can be
used in order to compute the Novikov - Shubin invariants.

In \cite{Fa3} we showed that the construction of the extended abelian category and extended
$L^2$ cohomology applies in many cases when one deals with {\it finite von Neumann categories}.
The category of Hilbert representations of a descrete group (considered in \cite{Fa1}, \cite{Fa2})
is a specific example of finite von Neumann category. Other examples, studied in \cite{Fa4},
include description of the geometry of growth and the behavior of the Betti numbers under
growth (L\"uck \cite{L4} type theorems).

The main purpose of the 
present paper is to classify Hermitian forms on torsion objects of a finite von Neumann category.
In the future publications we will apply the results obtained here to study the spectrum near zero of the
Dirac operator acting on the universal covering (the Novikov - Shubin signatures); we will
construct a combinatorial analogue of the eta-invariant; also, we will apply the technique of this paper
to study the knot concordance problem (in particular we will get a version of the Casson - Gordon
knot cobordism invariant).

The definition of a Hermitian form on torsion Hilbertian modules is not straightforward.
Indeed, such torsion module is not a set, only an object of certain category. In this paper
we adopt the general formalism of Hermitian forms in categories, as developed
by H.-G. Quebbemann, W. Scharlau, M. Schulte \cite{QSS}, A. Ranicki
\cite{R1, R2}, C.T.C. Wall \cite{W1, W2} and others. We review this formalism in \S 1. 
In \S 3
a duality in the category of torsion Hilbertian modules is described; using the duality, one may
apply the general formalism of Hermitian forms in this situation.

In section \S 4 we show that any non-degenerate torsion Hermitian form can be represented as a 
{\it discriminant form} of a {\it degenerate form} on a projective object, Theorem 4.7. 
Because of this, the study of non-degenerate torsion forms can be 
viewed as {\it a study of degenerate projective forms modulo the 
non-degenerate ones}. In this section we also describe completely the
equivalence relation (excision) between degenerate projective forms which corresponds 
to congruence of their discriminant forms, cf. Theorem 4.8. Hence, Theorems 4.7 and 4.8
give a classification of torsion Hermitian forms. 

A slightly similar approach to classification of degenerate Hermitian forms was suggested recently by
E. Bayer - Fluckiger and L. Fainsilber \cite{BF}. They showed that the classification of  
(degenerate) forms can be reduced to studying non-degenerate forms in the category of morphisms; 
note that the extended abelian category is a factor-category of the category of morphisms
(we factor out the null-homotopic morphisms). 

In \S 5 we introduce {\it positively and negatively definite} torsion Hermitian forms. It turns out that
any non-degenerate torsion Hermitian form is a direct sum of a positively and a negatively definite
forms (Theorem 5.3). It is much harder to establish the uniqueness of this decomposition; we prove
the uniqueness in \S 7 under an additional assumption that the initial von Neumann category is 
{\it superfinite} (Theorem 7.7). The definition of superfiniteness (cf. 7.2) is analogous to the 
von Neumann definition of finiteness with the only difference that here we care about the torsion 
objects instead of projective ones, as in the classical case. We show some examples when the 
superfiniteness holds; these examples include the categories which describe the growth, 
as in \cite{Fa4}. 
Another important property of superfinite categories is given in Theorem 7.4; it states that a positively
definite form is fully determined by its module.

We also study metabolic and hyperbolic torsion Hermitian forms. We show that any form is metabolic
and the metabolizer is sometimes unique. 
But, we show, that the hyperbolicity imposes very strong requirements on the structure (the positive
and the negative parts must be identical), cf. Theorem 7.10.

I would like to thank my colleagues I. Gohberg and V. Matsaev for a number of useful conversations.

\vskip 1cm
\centerline{\bf Contents}

\S 1. Dualities and Hermitian forms in categories.

\S 2. Hermitian forms in von Neumann categories.

\S 3. Duality for torsion Hilbertian modules.

\S 4. Classification of Hermitian forms on torsion Hilbertian modules.

\S 5. Positively and negatively definite torsion Hermitian forms.

\S 6. Metabolic and Hyperbolic forms.

\S 7. Superfinite von Neumann categories.

\heading{\bf \S 1. Dualities and Hermitian forms in categories}\endheading 

In this section we will briefly review
a general formalism of studying Hermitian forms in additive categories,
as developed in \cite{QSS, R1, R2}. We make some adjustments and modifications,
which will be convenient for the subsequent sections.
In particular, we emphasize the role of the group of automorphisms of 
a category with duality.

In the rest of this article we will apply this general formalism to study Hermitian forms
in von Neumann categories and on torsion objects of their extended abelian categories, 
cf. \S\S 2 - 7.

In all these applications we will always have $\C$ (the field of complex numbers) as part of the 
automorphism group of the category; therefore we will ignore the distinction 
between the quadratic and the symmetric cases as in \cite{R1, R2}.

\subheading{1.1. Dualities in categories}
Let $\CC$ be a category. We do not assume that $\CC$ is abelian or additive,
if not indicated explicitly.

\subheading{Definition} {\it A duality in category $\CC$ is a pair $(D, s)$, 
consisting of a contravariant 
functor $D: \CC\to \CC$ and an isomorphism
of functors $s: \Id \to D\circ D$, satisfying the following condition:
for any object $M$ of category $\CC$ the morphisms
$D(s_M)$ and $s_{D(M)}$ are inverse to each other:}
$$s_{D(M)}\ =\ D(s_M)^{-1}.\tag1-1$$

Here $s_M$ denotes the natural isomorphism $M\to D(D(M))$
determined by the isomorphism of functors $s$. In particular, $s_{D(M)}$ is a 
morphism $D(M)\to D(D(D(M)))$. The morphism 
$D(s_M):D(D(D(M)))\to D(M)$ is obtained
by applying the duality functor $D$ to $s_M$.

The object $D(M)$ is said to be {\it the dual} of $M$; for any morphism
$f:M\to N$ in $\CC$ we have {\it the dual morphism} $D(f):D(N)\to D(M)$.

\subheading{Remark} A. Ranicki \cite{R1, R2} calls duality 
$(D,s)$
{\it an involution} on the category $\CC$; he also assumes (as well as in \cite{QSS}) 
that there is given a canonical 
identification of any object $M$ with $D(D(M))$ via $s_M$, so that the isomorphism of functors
$s$ becomes the identity. In this paper we prefer not doing this.
In fact, we will deal here with situations when there are several
different dualities in different subcategories of a given abelian category, so mentioning both $D$ and 
$s$ is essential.

The following simple observation will be useful in the sequel.

\proclaim{Lemma} If $\CC$ is an abelian with duality $(D, s)$ so that the functor $D$ is additive.
Then $D$ is exact, i.e. it maps any exact sequence into an exact sequence.
\endproclaim
\demo{Proof} First, we observe that the map $D:\hom(M,N)\to 
\hom(D(N),D(N))$ is bijective. If $D$ is additive, then it maps zero morphism
to zero; thus we obtain that $f:M\to N$ is a monomorphism if and only if 
$D(f)$ is an epimorphism.
It follows that $D$ maps kernels to cokernels and conversely.
If $0\to N_1 @>\alpha>> N_2 @>\beta>> N_3\to 0$ is an exact sequence,
then $\beta$ is cokernel of $\alpha$ and therefore $D(\beta)$ is kernel
of $D(\alpha)$ which means that 
$$0\to D(N_3) @>D(\beta)>> D(N_2) @>D(\alpha)>> D(N_1)\to 0$$
is exact. \qed
\enddemo

\subheading{1.2. The symmetry group} Given 
a category $\CC$ with duality 
$(D,s)$, we will consider {\it its symmetry group} $G(\CC)$ , which is defined
as the set of all natural isomorphisms of the identity functor
$\Id_{\CC}$ to itself which are compatible with the duality (in the sense 
explained below). An element $\epsilon\in G(\CC)$ assigns an isomorphism 
$\epsilon_M:M\to M$
to any object $M$ of $\CC$ so that for any morphism $f:M\to N$ we have
$f\circ \epsilon_M\ =\ \epsilon_N\circ f$. We say that $\epsilon$ {\it is 
compatible with duality} if for any $M\in \ob(\CC)$ holds
$$s_{M}\circ \epsilon_{M}\ =\ DD(\epsilon_{M})\circ s_{M}.\tag1-2$$ 
In other words, we require that
the following diagram commutes
$$
\CD
M@>{\epsilon_M}>>M\\
@V{s_M}VV       @VV{s_M}V\\
DD(M)@>>{DD(\epsilon_M)}>DD(M).
\endCD\tag1-3
$$
This is equivalent to 
$$DD(\epsilon_M)\ =\ \epsilon_{DD(M)}.\tag1-4$$

The group $G(\CC)$ may clearly depend on the choice of duality $(D,s)$ in the category 
$\CC$ although we will not indicate this fact in our notation.

\subheading{Remark} It is easy to see that the group $G(\CC)$ of symmetries of any category 
with duality is always {\it abelian}.

There is naturally defined {\it involution in} $G(\CC)$: 
if $\epsilon\in G(\CC)$, let us denote by 
$\overline\epsilon\in G(\CC)$ the natural transformation 
$\Id_{\CC}\to \Id_{\CC}$
given by the diagram 
$$
\CD
M@>{\overline\epsilon_M}>>M\\
@V{s_M}VV       @VV{s_M}V\\
DD(M)@>>{D(\epsilon_{D(M)})}>DD(M)
\endCD.\tag1-5
$$
One easily checks that this diagram defines a natural transformation 
$\overline\epsilon: \Id_\CC\to \Id_\CC$ and thus an element of the symmetry
group $G(\CC)$. We observe the following useful formulae
$$
\aligned
\epsilon_{DD(M)} = &D(\overline\epsilon_{D(M)}),
\quad \overline{\overline\epsilon}\ = \epsilon, \,\\
\overline{\epsilon_1\circ\epsilon_2}\ = 
\ \overline{\epsilon_2}\circ 
\overline{\epsilon_1}\ = &\ 
\ \overline{\epsilon_1}\circ 
\overline{\epsilon_2}
\quad\text{for any}\quad \epsilon_1,\epsilon_2\in G(\CC).
\endaligned\tag1-6$$

\subheading{1.3. Hermitian forms} 
Given a category $\CC$ with duality $(D,s)$, one studies 
{\it Hermitian forms} in $\CC$. 

\subheading{Definition} {\it A form in a category with duality $\CC$
is a morphism
$\phi :M\to D(M)$, where $M$ is an object of $\CC$. $\phi$ is
called non-degenerate if it is an isomorphism in $\CC$.
Given a form $\phi:M\to D(M)$, consider the composition
$$\phi^\dagger\ \  :M@>{s(M)}>>D(D(M))@>{D(\phi)}>>D(M).\tag1-7$$
The form $\phi^\dagger$ is called the transpose of $\phi$.}

Note that 
$$\phi^{\dagger\dagger}\ =\ \phi.\tag1-8$$

A form $\phi:M\to D(M)$ is called $\epsilon$-{\it Hermitian}, where 
$\epsilon$ is an symmetry of the category with duality $\CC$, i.e. 
$\epsilon\in G(\CC)$, if 
$$\phi^\dagger\ =\ \epsilon_{D(M)}\circ \phi\ =\ 
\phi \circ \epsilon_{M}.$$
We will write this briefly as 
$\phi^\dagger \ =\ \epsilon\phi.$
One easily checks that for $\eta\in G(\CC)$, $(\eta\phi)^\dagger = \overline \eta
\phi^\dagger$. Hence if $\phi$ is $\epsilon$-Hermitian then 
$$\phi = \phi^{\dagger\dagger} = (\epsilon\phi^\dagger)^\dagger = 
\overline \epsilon \phi^\dagger = \epsilon\overline\epsilon \phi.$$

\subheading{Conclusion} Non-degenerate $\epsilon$-Hermitian forms exit only if the symmetry $\epsilon\in G(\CC)$ 
is unitary:
$\epsilon\overline\epsilon\ =\ 1.$

Let $\phi_1:M_1\to D(M_1)$, and 
$\phi_2:M_2\to D(M_2)$ be two forms in category with duality $\CC$. 
They are
called {\it congruent (or isometric)} if there is an isomorphism 
$f:M_1\to M_2$ such
that $\phi_1\ =\ D(f)\circ\phi_2\circ f$.

Congruence is an equivalence relation.
If two forms are congruent and one of them is $\epsilon$-Hermitian for 
some $\epsilon\in G(\CC)$, then the other is also $\epsilon$-Hermitian.

Now we will show how the problem of describing the congruence classes of $\epsilon$-Hermitian
forms may depend on $\epsilon$.
Suppose that $\eta\in G(\CC)$ and $\phi:M\to D(M)$ is an $\epsilon$-Hermitian  
form for some $\epsilon\in G(\CC)$. Consider the form $\psi=\eta\phi:M\to D(M)$.
Then an easy computation shows that $\psi$ is 
$\epsilon\overline\eta\eta^{-1}$-Hermitian. It follows from the above remark 
that {\it the set of congruence classes of $\epsilon$-Hermitian
forms in $\CC$ depends only on the class of $\epsilon$ in
the factor-group}
$$\E(\CC)\ = \ 
\{\epsilon\in G(\CC); \epsilon\overline\epsilon=1\}/1\}/\{\overline\eta\eta^{-1};
\eta\in G(\CC)\}.\tag1-9$$
We will call (1-9) {\it the group of types of forms} in $\CC$.

Note that any element in $\E(\CC)$ is of order two:  
$\epsilon\overline\epsilon=1$ implies that 
$\epsilon^2\ =\ \overline\eta\eta^{-1}$ for $\eta\ =\ \overline\epsilon$.

\subheading {1.4. Example: Duality of finitely generated projective modules}
This is the most familiar example of a category with duality.

Let $\Lambda$ denote a ring with involution which will be denoted 
$\lambda\mapsto \overline\lambda$.
Consider the category $\lm$ of left finitely generated 
projective $\Lambda$-modules and $\Lambda$-homomorphisms.
If $M$ is a projective finitely generated $\Lambda$-module, then by $D(M)$
we will denote the set of all {\it anti-linear} homomorphisms 
$f:M\to \Lambda$, i.e. all additive maps satisfying 
$$f(\lambda m)\ =\ f(m)\overline\lambda\tag1-10$$
for all $m\in M$ and $\lambda\in\Lambda$. 
A left $\Lambda$-module structure on $D(M)$ is defined as follows:
for $\lambda\in\Lambda$ and $f\in D(M)$, 
$$(\lambda\cdot f)(m)\ =\ \lambda f(m)\tag1-11$$
where $m\in M$. It is easy to see that this formula
defines a structure of left $\Lambda$-module on $D(M)$. One checks also
that $D(M)$ is projective and finitely generated.
Any $\Lambda$-homomorphism $\phi:M\to N$ between finitely generated
projective left $\Lambda$-modules induces naturally a homomorphism
$$D(\phi):D(N)\to D(M)$$
by the rule $D(\phi)(f)=f\circ \phi$, where $f\in D(N)$. 
Thus we obtain a contravariant functor $D: \lm \to \lm$. 

The canonical isomorphism
$$s_M: M \ \to\  DD(M)\tag1-12$$
is given by 
$$m\mapsto (f\mapsto \overline{f(m)}),\quad f\in D(M).\tag1-13$$
One checks that, firstly, $(f\mapsto \overline{f(m)})$ is an anti-liner 
functional on $D(M)$, i.e. an element of $DD(M)$, and, secondly, the map
(1-12) is a $\Lambda$-isomorphism. The pair $(D,s)$ is a duality in category
$\lm$.

The group $G(\CC)$ of symmetries of this category with duality coincides 
with the set of all invertible elements in the center $Z(\Lambda)$ of the 
ring $\Lambda$, i.e. $G(\CC)\ =\ Z(\Lambda)^\ast.$ 
The involution on $G(\CC)$, defined 
as in subsection 1.2 (cf. (1-5)), coincides with the
restriction of the original involution of $\Lambda$ onto $Z(\Lambda)^\ast$. 
We obtain that the group of types of forms is
$$\E(\CC)\ = \ 
\{\epsilon\in Z(\Lambda); \epsilon\overline\epsilon=1\}/
\{\overline\eta\eta^{-1};\eta\in Z(\Lambda)^\ast\}.\tag1-14$$
In particular, we obtain that {\it the group $\E(\CC)$ depends
only on the center $Z(\Lambda)$ and the action of the involution of 
$\Lambda$ on
the center.}

\heading{\bf \S 2. Hermitian forms in von Neumann 
categories}\endheading

In this section we will first recall the main properties of von Neumann categories,
which were introduced in \cite{GLR} and independently in \cite{Fa3}. Our exposition
will follow \cite{Fa3}, \S 2, \S 5. 
We show that the classification of non-degenerate Hermitian forms
in a von Neumann category is immediate. 

\subheading{2.1. Hilbertian spaces} Recall that {\it a Hilbertian space} 
(cf. \cite{P}) is a topological vector space $H$, which
is isomorphic to a Hilbert space in the category of topological vector 
spaces. In other words, there exists a scalar product on $H$, such that 
$H$ with this scalar product is a Hilbert space with the originally given 
topology. Such scalar
products are called {\it admissible}. Given one admissible scalar product
$\< \ ,\ \>$ on $H$, any other admissible scalar product is given by
$$\<x,y\>_1\ =\ \<Ax,y\>,\quad x,y\in H,\tag2-1$$
where $A: H\to H$ is an invertible positive operator $A^\ast=A,\quad A>0$.
Hilbertian spaces naturally appear as Sobolev spaces of sections of
vector bundles,
cf. \cite{P}.

Let us denote by $\Hilb$ the category of Hilbertian
spaces and continuous linear maps. This category has an obvious duality $(D,s)$, where
$D(H)$ is defined as the space $H^\ast$ of all 
anti-linear continuous functionals on $H$ (i.e the set of all 
continuous $\R$-linear maps $\phi:H\to \C$, such that 
$\phi(\lambda h)=\overline \lambda\phi(h)$ for all $\lambda\in \C$ and 
$h\in H$; here the bar denotes the complex conjugation).
We consider the action of $\C$ on 
$D(H)=H^\ast$ given by
$(\lambda\cdot \phi)(h)=\phi(\overline\lambda^\ast\cdot h)$ for all $h\in H$. 
The canonical isomorphism 
$$s_{H}: H\to DD(H)=H^{\ast\ast}\tag2-2$$ 
is $h\mapsto(\phi\mapsto\overline{\phi(h)})$, where $h\in H$, and 
$\phi\in H^\ast$. One checks that the condition (1-1) is satisfied and 
so we have a category with duality.

\subheading{2.2. Hilbertian von Neumann categories} 
{\it A Hilbertian von Neumann category} is an additive subcategory $\CC$
of $\Hilb$ with the following properties: \roster
\item {\it for any $H\in \ob(\CC)$ the dual space $H^\ast$ is also an object of $\CC$ and there
is a $\CC$-isomorphism $\phi: H\to H^\ast$ such that the formula
$$\<x,y\> = \phi(x)(y), \quad x, y\in H$$
defines an admissible scalar product on $H$;
\item for any $H\in \ob(\CC)$ the isomorphism (2-2) also belongs to $\CC$;
\item the adjoint of any morphism in $\CC$ also belongs to $\CC$;
\item the kernel $\ker f=\{x\in H; f(x)=0\}$ 
of any morphism $f: H\to H'$ in $\CC$ and the natural inclusion $\ker f \to H$ belong 
to $\CC$;
\item for any $H, H'\in \ob(\CC)$ the set of morphisms
$$\hom_{\CC}(H, H')\subset \hom_{\Hilb}(H, H') = \L(H, H')$$
is a linear subspace closed with respect to the weak topology.}
\endroster
 
Hence, objects of $\Cal C$ have structure of Hilbertian spaces and possibly some additional
structure,
and morphisms of $\CC$ are (faithfully) represented by bounded linear maps. 

Condition (5) is similar to the well-known condition in the definition of von Neumann algebras,
which explains our term.
Recall, that the weak topology on the space of bounded linear operators
$f:H\to H'$ is given by the family of seminorms 
$$p_{\phi,x}(f)\ =\ |\<\phi,f(x)\>|,\quad \text{where}\quad \phi\in H', 
\quad x\in H.\tag2-3$$

Note also, that given any object $H\in\ob(\CC)$, a choice of $\CC$-admissible
scalar product on $H$, determines an involution on the algebra 
$\hom_{\CC}(H,H)$ (adjoint operator) and the space $\hom_{\CC}(H,H)$ 
considered with this involution is a von Neumann algebra.

The above conditions (1) - (5) imply the following properties:
\roster
\item"{(6)}" {\it The closure of the image $\cl(\im f)$ of any morphism 
$f: H\to H'$ in $\CC$ and also the natural projection 
$H'\to H'/\cl(\im f)$ belong to $\CC$.
\item"{(7)}" Suppose that $H'\subset H$ is a closed subspace. 
If $H'$, $H$ and the inclusion $H'\to H$ belong to $\CC$ then the orthogonal
complement ${H'}^\perp$ with respect to an admissible scalar product on
$H$ and the inclusion ${H'}^\perp\to H$ belong to
$\CC$.}
\endroster

\subheading{Remark} The duality structure in $\Hilb$, described in 2.1, induces a duality structure
in any von Neumann category $\CC \subset \Hilb$. This follows from conditions (1), (2) and (3)
of subsection 2.2. Hence {\it any Hilbertian von Neumann category is a category with duality.}

\subheading{2.3. Example} 
The simplest example of a Hilbertian von Neumann category is the following.
Let $\A$ be an algebra over 
$\C$ with involution, which on $\C$ coincides with the complex conjugation. 
{\it A Hilbertian representation of $\A$ } is a Hilbertian 
topological vector space $H$
supplied with a left action of $\A$ by continuous linear maps 
$\A\to \L(H,H)$. 
A morphism $f: H\to H'$ between two Hilbertian representations of $\A$ is 
defined as a bounded linear map commuting 
with the action of the algebra $\A$.
There is a canonical duality in the category 
of all Hilbertian representations of a given $\ast$-algebra $\A$.
Namely, given a Hilbertian representation $H$, consider the space
$D(H)=H^\ast$ 
of all anti-linear continuous functionals on $H$. 
Consider the following action of $\A$ on 
$D(H)=H^\ast$: if $\phi\in H^\ast$ and $\lambda\in \A$ then 
$(\lambda\cdot \phi)(h)=\phi(\lambda^\ast\cdot h)$ for all $h\in H$. 
The canonical isomorphism 
$$s_{H}: H\to DD(H)=H^{\ast\ast}\tag2-4$$ 
is given by 
$h\mapsto(\phi\mapsto\overline{\phi(h)})$, where $h\in H$, and 
$\phi\in H^\ast$. Category of all Hilbertian representations of $\A$ is a von Neumann category.

\subheading{2.4. Finite von Neumann categories}
An object $H$ of a Hilbertian von Neumann category $\CC$ will be called {\it 
finite} if the only closed $\CC$-submodule $H_1\subset H$ which is 
isomorphic to $H$ in $\CC$ is $H_1=H$.

A Hilbertian von Neumann category will be called {\it finite} iff all its  
objects are finite.

We will use the following property of finite categories.

\proclaim{2.5. Proposition} Let $\CC$ be a finite von Neumann category and let
$\alpha: H \to H'$ be an injective morphism of $\CC$
with dense image. Then $H$ is isomorphic to $H'$ in $\CC$. For any injective nonzero
morphism $\beta: H'' \to H'$ the preimage $\beta^{-1}(\alpha(H)) \subset H''$
is dense (and so non-empty).
\endproclaim
\demo{Proof} The first statement was proven in \cite{Fa1}, Lemma 2.3. By Lemma 2.4 of
\cite{Fa1} the preimage $\beta^{-1}(\alpha(H)) \subset H''$ is not empty. If it is not
dense, consider a complement $H''' \subset H''$ to the closure of $\beta^{-1}(\alpha(H))$.
Let $\beta' : H''' \to H'$ be the restriction of $\beta$. Then the preimage
${\beta'}^{-1}(\alpha(H))$ is empty, in contradiction with the established above
property. \qed
\enddemo

\subheading{2.6. More examples of von Neumann categories}.

{\it Example 1.} The following is one of the most important examples of 
von Neumann categories. Let $\B$ be a von Neumann algebra acting on a 
Hilbert space $H$. 
Denote by $\A$ be the commutant of $\B$. We will consider the following 
category $\CC$ of Hilbertian representations of $\A$. Objects of $\CC$  
are in one-to-one correspondence with projections $e\in M(n)\otimes \B$, 
$e^2=e$, for some $n$, 
where $M(n)$ is the $n\times n$-matrix algebra.
For each projection $e$ the corresponding Hilbert representation of
$\A$ is $e(\C^n\otimes H)$. If $e_1$ and $e_2$ are two projections, 
then the set of morphisms 
$\hom_{\CC}(e_1(\C^{n_1}\otimes H),e_2(\C^{n_2}\otimes H))$ 
is the set of all bounded linear maps of the form
$e_2be_1:e_1(\C^{n_1}\otimes H)\to e_2(\C^{n_2}\otimes H)$ 
where $b$ is given by an $n_1\times n_2$-matrix with entries in $\B$. 
$CC$ is a von Neumann category.

{\it Example 2.} Let $Z$ be a locally compact Hausdorff space and let $\mu$ be a 
positive Radon measure on $Z$. Let $\A$ denote the algebra 
$L^\infty_{\C}(Z,\mu)$ (the space of essentially bounded $\mu$-measurable
complex valued functions on $Z$, in which two functions equal locally
almost everywhere, are identical). The involution on $\A$ is given by the
complex conjugation. 
We will construct a category $\CC$ of Hilbert representations
of $\A$ as follows. The objects of $\CC$ are in one-to-one correspondence
with the $\mu$-measurable fields of 
finite-dimensional Hilbert spaces $\xi\to\H(\xi)$ over $(Z,\mu)$, 
cf. \cite{Di}, part II, chapter 1. For any such measurable field of Hilbert
spaces, the corresponding Hilbert space is the direct integral
$$H = \int^{\oplus}\H(\xi)d\mu(\xi)\tag2-5$$
defined as in \cite{Di}, part II, chapter 1. The algebra $\A$ acts on the 
Hilbert space (2-5) by pointwise multiplication. 

Suppose that we have two $\mu$-measurable finite-dimensional fields of
Hilbert spaces $\xi\to\H(\xi)$ and $\xi\to\H'(\xi)$ over $Z$. Then we 
have two corresponding Hilbert spaces, $H$ and $H'$, given as direct
integrals (2-5). We define the set of morphisms $\hom_{\CC}(H,H')$ as the
set of all bounded linear maps $H\to H'$ given by {\it decomposable linear
maps} 
$$T = \int^\oplus T(\xi)d\mu(\xi),\tag2-6$$
where $T(\xi)$ is {\it an essentially bounded measurable field of linear 
maps}
$T(\xi):\H(\xi)\to\H'(\xi)$, cf. \cite{Di}, part II, chapter 2.  

The kernel of any decomposable linear map as above can be represented as the
direct integral of a finite-dimensional field of Hilbert spaces and so
condition (c) of section 2.4 is satisfied. Condition (a) of section
2.4 is also satisfied since the adjoint of the map $T$ given by (2-6) is
$$T^\ast=\int^\oplus T(\xi)^\ast d\mu(\xi)\tag2-7$$
(by \cite{Di}, part II, chapter 2, \S 3, Proposition 3).
Condition (d) of section 2.4 is satisfied as follows from Theorem 1
of \cite{Di}, part II, chapter 2, \S 5.

Thus, we obtain a von Neumann category.
This category is finite.
Note that this category $\CC$ depends only on the {\it class of the measure
$\mu$}.

Other examples of von Neumann categories can be found in \cite{Fa3}, \S 2.

\subheading{2.7. Group of types of forms}
Since any von Neumann category $\CC$ is a category with duality, we may consider 
its symmetry group $G(\CC)$, cf. 1.2. The corresponding group of types
of forms $\E(\CC)$ (cf. (1-9)) is trivial in all examples mentioned above.
For instance, in example 1 of section 2.6 the symmetry group $G(\CC)$
coincides with the group $Z(\A)^\ast$ of invertible elements in the center
of the von Neumann algebra $\A$. If $\A$ is a factor, then $Z(\A)=\C$ and
$\E(\CC)=0$. If $\A$ is not a factor, then $Z(\A)$ can be represented
as the ring $L^\infty_{\C}(Z,\mu)$ of essentially bounded measurable functions
on a measure space $(Z,\mu)$. Then any unitary element $\epsilon\in Z(\A)$
is represented by a function with values in the circle $S^1$. Given such 
$\epsilon$, it is clear that one can find another function $\eta$ on $Z$
with values in $S^1$ such that $\eta^2 = \epsilon$. Thus $\E(\CC)=0$.

\subheading{2.8. Hermitian forms in von Neumann categories} Let $\CC$ be a
finite von Neumann category. It is a category with duality (cf. Remark in 2.2)
and so we may apply to $\CC$ the formalism developed in \S 1. 

Any Hermitian form $\phi: H\to H^\ast$ in $\CC$ determines a 
{\it continuous scalar product}
$$\<\ ,\ \>_\phi: H\times H\to \C,$$
given by the formula
$$\<h,h^\prime\>_\phi \ = \ \phi(h)(h^\prime),\quad\text{for}\quad 
h, h^\prime\in H,$$
which satisfies the following properties:
$$\align
\<h,h^\prime\>_\phi \ &=\ \overline{\<h^\prime,h\>}_\phi\qquad
\text{(because $\phi$ is Hermitian)},\\
\<\lambda h,h^\prime\>_\phi\ &= \ \<h,\lambda^\ast h^\prime\>_\phi,
\quad \lambda\in\A, \quad\text{(since $\phi$ commutes with $\A$),}\tag2-8\\
\<\lambda h,h^\prime\>_\phi \ &=\ \lambda\<h,h^\prime\>_\phi\quad\text{for}
\quad\lambda\in \C.\\
\<h,h'\>_\phi &\le C\cdot ||h||\cdot ||h'||,\quad\text{(since $\phi$ is continuous).}
\endalign
$$
Here we assume that we have chosen a $\CC$-admissible scalar product $\<\ ,\ \>$ on $H$
(i.e. an admissible scalar product which appears in 2.2.(1))
and that $||\cdot||$ denotes the corresponding norm.

Note that if the kernel $\phi^{-1}(0)$ is zero, then the image 
$\phi(H) \subset H^\ast$ is dense, as follows from finiteness of $H$. 

If $\phi$ is Hermitian then $\<h,h\>_\phi$ is real.
We will say that a form $\phi$ as  above is {\it positively definite} if 
$\<h,h\>_\phi$ is a positive real number for all $h\in H, \ h\ne 0$. 
Similarly, we will say that a form $\phi$ is
{\it negatively definite} if $-\phi$ is positively definite.

We will call a Hermitian form $\phi: H \to H$ {\it weakly invertible}
is it is injective $\phi^{-1}(0) = 0$. 

\proclaim{2.9. Proposition} In a finite von Neumann category $\CC$, any 
weakly non-degenerate Hermitian form $\phi: H\to H^\ast$ 
can be represented as an orthogonal sum 
$$\phi\ =\ \phi_+\oplus\phi_-,\quad H=H_+ \oplus H_-\tag2-9$$
of a positively definite form $\phi_+: H_+\to H_+^\ast$ and a negatively 
definite form $\phi_-: H_-\to H_-^\ast$. 
\endproclaim

\demo{Proof} Fix a $\CC$-admissible scalar product $\<\ ,\ \>$ on $H$. Then 
the form $\phi$ determines a bounded self-adjoint injective operator $A_\phi\in
\hom_{\CC}(H,H)$ as above. By the spectral theorem 
we have the following decomposition
$$A_\phi\ =\ \int_{-\infty}^\infty \lambda dE_\lambda,$$
where $E_\lambda$ are the spectral projectors determined by $A_\phi$. Note
that all these projectors $E_\lambda$ belong to $\hom_{\CC}(H,H)$; in
particular, they commute with the action of $\A$.
Consider the mutually orthogonal projections
$P_- = E_0$ and $P_+ = E_N - E_0,$
where $N$ is sufficiently large, and the corresponding orthogonal
decomposition
$$H\ =\ H_+\oplus H_-,\quad\text{where}\quad H_\pm = P_\pm H.\tag2-10$$
The restriction of the form $\phi$ 
on $H_\pm$ is positively (negatively) definite. \qed

\enddemo

\proclaim{2.10. Proposition} Suppose that a Hermitian form $\phi: H\to H^\ast$ 
in a finite von Neumann category $\CC$ is non-degenerate (i.e. $\phi$ is an isomorphism
in $\CC$). Then the decomposition
(2-9) is unique. This means that for any other representation of $H$ as the direct sum of
two closed subspaces $H=H'_+ \oplus H'_-$, such that $H'_\pm\in\ob(\CC)$, and
the form $\phi$ is positively definite on $H'_+$
and negatively definite on $H'_-$, and $H'_+$ is orthogonal to $H'_-$ with respect to
$\phi$, then $H_+$ is isomorphic to $H'_+$ and $H_-$ is
isomorphic to $H'_-$ in $\CC$; moreover, the restriction $\phi|_{H_\pm}$ is congruent to
$\phi|_{H'_\pm}$. 
\endproclaim
\demo{Proof} Let $L$ be an arbitrary closed subspace of $H$ such that the form $\phi$ is
positively definite on $L$. Each $x\in L$ can be uniquely represented as $x = x_+ + x_-$ with
$x_\pm \in H_\pm$. Consider the projection $\pi : L \to H_+$, where $x\mapsto x_+$.
We claim that this map $\pi: L\to H_+$ is an injection with closed image.
In fact, if there exists
a sequence $(x_n)$ with $x_n\in L$, $||x_n|| = 1$, and $\pi(x_n)\to 0$, then
we obtain
$$0 < \ \<x_n,x_n\>_\phi \le \<\pi(x_n),\pi(x_n)\>_\phi \to 0,$$
and therefore $\<x_n,x_n\>_\phi \to 0$,
which contradicts the assumption that $\phi$ is non-degenerate.

Thus, if we have two decompositions $H=H_+ \oplus H_-$ and $H=H'_+ \oplus H'_-$,
then $H_+$ can be mapped into $H'_+$ and $H'_+$ can be mapped into $H_+$. From finiteness
of $\CC$ it now follows, that $H_+$ is isomorphic to $H'_+$ in $\CC$.

To show that $\phi|_{H_\pm}$ is congruent to $\phi|_{H'_\pm}$
we may use an easy fact that {\it any two non-degenerate positively definite forms on a given
object of $\CC$ are congruent}. \qed
\enddemo

\heading{\bf \S 3. Duality for torsion Hilbertian modules}\endheading

In this section we construct a duality in the abelian category of torsion Hilbertian modules.

\subheading{3.1. Abelian extension of a von Neumann category} Given a Hilbertian
von Neumann category $\CC$, 
there is a bigger category $\E(\CC)$, containing $\CC$ 
as a full subcategory \cite{Fa}, \cite{Fa1}. 
The advantage of $\E(\CC)$ is that it is {\it an abelian category.}
A brief description of the construction is given 
below; we refer to  \cite{Fa3} for more details. 

{\it An object} of the category $\E(\CC)$ is defined as a
morphism $(\alpha:A^\prime\to A)$ in the original category $\cc$. 
Given a pair of objects $\X=(\alpha: A^\prime \to A)$ and 
$\Y=(\beta:B^\prime\to B)$ of $\E(\cc)$,
a {\it morphism} $\X\to\Y$ in category $\E(\cc)$ is an  
equivalence class of morphisms $f:A\to B$ of category $\CC$
such that $f\circ\alpha=\beta\circ g$
for some morphism $g:A^\prime \to B^\prime$ in $\CC$. 
Two morphisms $f:A\to B$ and $f^\prime:A\to B$ of $\CC$ represent 
{\it identical morphisms $\X\to\Y$ of $\E(\CC)$} iff 
$f-f^\prime = \beta\circ F$ for some morphism $F:A\to B^\prime$ of category
$\CC$. 
The morphism $\X\to\Y$, represented by $f:A\to B$,
is denoted by
$$[f]:(\alpha:A^\prime\to A)\ \to\ (\beta:B^\prime\to B)\quad
\text{or by}\quad
[f]:\X\to\Y.\tag3-1$$
{\it Composition} of morphisms in $\E(\CC)$ is defined as composition of the 
corresponding morphisms $f$ in the category $\CC$.

The category $\E(\CC)$ is an abelian category, cf. \cite{Fa3}, Proposition 1.7. It is
called {\it the abelian extension of 
the category} $\CC$.

It is shown in \cite{Fa3}, section 1.4, that any object $\X$ of $\eca$
is isomorphic in $\eca$ to an object $(\alpha:A'\to A)$, where the 
morphism $\alpha$ is injective.

\subheading{3.2. Torsion Hilbertian modules} An object $\X=(\alpha:A'\to A)$ 
of the extended category $\eca$ is called {\it torsion} iff the image of 
$\alpha$ is dense in $A$. 

We will denote by $\tca$ the full subcategory of $\eca$ generated by all
torsion objects. $\tca$ is called {\it the torsion subcategory of $\eca$}.
It is shown in \cite{Fa3}, \S 3, that: 

{\it if $\CC$ is a finite von Neumann 
category, then $\tca$ is an abelian category.}

Objects of $\tca$ are called {\it torsion Hilbertian modules}.

\subheading{3.3. Duality for torsion Hilbertian modules} Our aim now is to 
construct duality in the category $\tca$ of torsion Hilbertian modules, i.e.
a contravariant functor 
$$\e:\tca \to \tca \tag3-2$$ 
and an isomorphism
of functors 
$$s:\Id\to \e\circ\e,\tag3-3$$ 
cf. 1.1.
For any torsion Hilbertian module $\X= (\alpha:A^\prime\to A)$ (where $\alpha$
is injective) define the {\it dual module} $\e(\X)$ by
$$\e(\X)\ =\ (D(\alpha):D(A)\to D(A^\prime)),\tag3-4$$
where $D$ denotes the duality functor in $\CC$,
cf. sections 2.1, 2.2.

Suppose now that we have two torsion objects $\X=(\alpha: A^\prime\to A)$
and $\Y=(\beta: B^\prime\to B)$ with injective $\alpha$ and $\beta$ and let
$[f]:\X\to\Y$ be a morphism represented by a diagram
$$
\CD
(A^\prime@>\alpha>>A)\\
 & &                @VVfV\\
(B^\prime@>\beta>>B).
\endCD\tag3-5
$$
According to definition \S 3.1, there exists a morphism 
$h:A^\prime\to B^\prime$
making this diagram commutative; this $h$ is in fact unique, because of
injectivity of $\beta$. We define the {\it dual morphism of 
$[f]$} as the morphism
$$\e([f])\ =\ [D(h)] : \e(\Y)\ \to \e(\X).\tag3-6$$
It is represented by the diagram
$$
\CD
(D(B)@>{D(\beta)}>>D(B^{\prime}))\\
& &      @VV{D(h)}V\\
(D(A)@>>{D(\alpha)}>D(A^{\prime})).
\endCD\tag3-7
$$

Since different morphisms $f$ may represent the same morphism in $\tca$,
cf. \S 3.1, we have to check correctness of the above definition. 
If $F:A\to B^\prime$
is an arbitrary morphism, then the morphism $f^\prime=f+\beta\circ F$
represents the same morphism $[f]=[f^\prime]$ in $\E(\CC)$. 
Then the corresponding 
morphism $h^\prime$ is $h^\prime= h+F\circ\alpha$ and thus 
$$D(h^{\prime})\ =\ D(h)\ +\ D(\alpha)\circ D(F),$$
which means that $D(h)$ and $D(h^{\prime})$ represent the same morphism in
$\tca$, cf. \cite{Fa3}.

The isomorphism of functors $s:\Id\to \e\circ\e$ is given as follows.
For any torsion Hilbertian module $\X=(\alpha: A^\prime\to A)$ define
$$s_{\X}:\X\to \e(\e(\X))\quad\text{by}\quad 
s_{\X}=[s_A],\tag3-8$$ 
where $s_A$ is the canonical
isomorphism $A\to DD(A)$ in the category $\CC)$, cf. sections 2.1, 2.2. We
have the following commutative diagram
$$
\CD
A^\prime@>\alpha>>A\\
@V{s_{A^\prime}}VV   @VV{s_A}V\\
DD(A^\prime)@>>{DD(\alpha)}>DD(A),
\endCD\tag3-9
$$
representing the morphism $s_{\X}:\X\to \e(\e(\X))$.

\heading{\bf \S 4. Classification of Hermitian forms on torsion Hilbertian 
modules}\endheading 

This section plays a central role in the paper; here we give a classification of non-degenerate 
Hermitian forms on torsion Hilbertian modules. 

We refer to the paper of E. Bayer - Fluckiger and L. Fainsilber \cite{BF} where a problem
of classification of degenerate Hermitian forms was studied. Their results are similar in the spirit
to the results of this section.

\subheading{4.1} Let $\CC$ be a finite von Neumann category and let $\tca$ 
be the torsion subcategory of the extended abelian category $\eca$, cf. 3.2. 
We will study Hermitian forms in $\tca$, considered with 
the duality described in \S 3.

Let $\X=(\alpha:A^\prime\to A)$ be a torsion object of $\E(\A)$ represented
by an {\it injective} morphism $\alpha$ with dense image. Then a 
Hermitian form on $\X$ is given by a morphism 
$$\phi:\X\to\e(\X),\quad\text{with}\quad \phi^\dagger = \phi\tag4-1$$
(cf. \S 1). According to the definitions of \S 3 and \cite{Fa3},
the morphism $\phi$
is represented by a commutative diagram
$$
\CD
A^\prime@>\alpha>>A\\
@V{h}VV     @VV{f}V\\
D(A)@>>{D(\alpha)}>D(A^\prime)
\endCD\tag4-2
$$
in $\CC$. Clearly, $f$ determines $h$ uniquely  (since $\alpha$ is injective and has dense image and
hence $D(\alpha)$ is also injective). The condition 
$\phi^\dagger = \phi$ means (translating the general definitions of 
\S 1) that 
$$f -  h^\dagger\ =\ D(\alpha)\circ F,\tag4-3$$
where $F:A\to D(A)$ is a morphism of $\CC$ and $h^\dagger$ denotes the 
composition
$$A@>{s_A}>>DD(A)@>{D(h)}>>D(A^\prime),$$
i.e. the transpose of $h$ in $\CC$, cf. (1-7). Note that (4-3) implies that
$F$ is anti-symmetric, i.e. $F^\dagger = - F$.

\subheading{Definition} Pair of morphisms $(f,h)$ as in diagram (4-2) above,
will be called {\it symmetric 
presentation of the Hermitian form $\phi:\X\to \e(\X)$} if 
$$f\ =\ h^\dagger.\tag4-4$$

\proclaim{4.2. Lemma} Any Hermitian form in $\tca$ admits a 
symmetric presentation (not unique). 
\endproclaim
\demo{Proof} Let a pair $(f,h)$ be an arbitrary presentation of a morphism
$\phi$ as above. Denote
$$f_1\ =\ f - \frac{1}{2}D(\alpha)\circ F\quad\text{and}\quad
h_1\ =\ h - \frac{1}{2} F\circ \alpha,$$
where $F$ satisfies (4-3).
Then the obtained pair of morphisms $(f_1,h_1)$ also represents $\phi$ 
and we have $f_1\ =\ h^\dagger_1$. \qed
\enddemo

\subheading{4.3} A symmetric presentation $(f,h)$ of a form 
$\phi:\X\to \e(\X)$ can be also viewed as follows. The map 
$f:A\to D(A^\prime)$ defines a pairing
$\<\ ,\ \>: A\times A^\prime\to \C$
(where $\<a,a^\prime\>\ =\ f(a)(a^\prime)$ for $a\in A, a^\prime\in A^\prime$)
which is continuous as a function of two variables, $\C$-linear with 
respect to the first variable and skew-linear with respect to the second
variable, and satisfies
$\<\lambda a,a^\prime\>\ =\ \<a,\lambda^\ast a^\prime\>$
for any $a\in A$, $a^\prime\in A^\prime$ and $\lambda\in \A$. 
Besides, symmetricity
implies that the induced form 
$\{\ ,\ \}:A^\prime\times A^\prime\to \C,$
where $\{x,y\}\ =\ \<\alpha(x),y\>$ for $x,y\in A^\prime$, is 
{\it Hermitian}, i.e. $\{x,y\}\ = \overline{\{y,x\}}$.

\subheading{4.4. Discriminant forms} The easiest
way of constructing Hermitian forms on torsion objects in 
$\tca$ consists in the following. Suppose that we are given a Hilbertian 
module $A\in\ob(\CC)$ and a Hermitian form in category $\CC$ (cf. \S 2)
$\alpha:A\to D(A).$
We will suppose that $\alpha$ is injective, but not necessarily
surjective. Then $\alpha$ represents a torsion object 
$$\X =\ (\alpha:A\to D(A))\in\ob(\tca)\tag4-5$$
and, moreover, we have the Hermitian form 
$$\phi_\alpha: \X\to \e(\X),\tag4-6$$
given by the following diagram
$$
\CD
A@>\alpha>>D(A)\\
@V{1}VV   @VV1V\\
A@>>{\alpha^\dagger}>D(A).
\endCD\tag4-7
$$
The form $\phi_\alpha$ is clearly non-degenerate.

Thus, {\it any Hermitian form $\alpha$ in category $\CC$ produces a 
non-degenerate Hermitian form on a torsion object in $\tca$}. We will call
the obtained form $\phi_\alpha$ the {\it discriminant form of} $\alpha$,
because of similarity of this construction with the well-known 
construction of discriminant forms in number theory, cf. \cite{N}.

It is clear that the discriminant form $\phi_\alpha$ of any non-degenerate 
form $\alpha$ is trivial (i.e. it has $\X=0$). 
Roughly, {\it the discriminant form $\phi_\alpha$ 
measures the "way of degeneration" of the form $\alpha$}.

The construction of discriminant form can be made even simpler if we suppose
that a $\CC$-admissible scalar product on a Hilbertian module
$A$ is specified. In this situation we may use this scalar product to 
identify the dual module $D(A)$ with
$A$ and so the Hermitian form $\alpha$ is now represented by a self-adjoint
operator 
$$\alpha:A\to A,\qquad \alpha^\ast\ =\ \alpha,\qquad \alpha\in\hom_{\CC}(A,A).\tag4-8$$
Here we have
$\X\ =\ (\alpha:A\to A)$ and the discriminant form $\phi_\alpha:\X\to\e(\X)$
will be represented by the commutative diagram
$$
\CD
A@>\alpha>>A\\
@V{1}VV   @VV1V\\
A@>>{\alpha}>A.
\endCD\tag4-9
$$

\subheading{4.5. Example} Consider the von Neumann algebra $\A=L^\infty(S^1)$
acting on $L^2(S^1)$ by multiplication. Let $\CC$ denote the von Neumann 
category, which is obtained by applying the construction of Example 2 in 
section 2.6 to the circle $S^1$ with the Lebesgue measure $\mu$. 

The space $L^2(S^1)$ has a canonical $\CC$-admissible scalar product. 
Therefore, any real valued function $\alpha\in L^\infty(S^1)$, which does 
not vanish a subset of positive measure, determines (by multiplication) an
operator
$\alpha: L^2(S^1)\to L^2(S^1),$
which gives (as explained in (4-9)) a Hermitian form $\phi_\alpha:\X \to \e(\X)$
in category $\tca$, where $\X\ =\ ( \alpha: L^2(S^1)\to L^2(S^1))$.

Note that the same form $\phi_\alpha\ =\ \phi_\beta$ in $\tca$
will be described by any other function $\beta\in L^\infty(S^1)$ of the form
$$\beta(z)\ =\ \alpha(z) \ +\ \alpha^2(z) F(z),\quad z\in S^1,\tag4-10$$
for arbitrary real valued $F\in L^\infty(S^1)$. 

Now we face the main question:
what do the functions $\alpha(z)$ and $\beta(z)$ have in common? In fact, choosing
$F(z)$ we may make $\beta(z)$ to be arbitrary {\it far from the zeros of} 
$\alpha(z)$.

To understand the situation, suppose for simplicity, that the real valued 
function $\alpha$ on $S^1$
is continuous and vanishes at a single point $z_0\in S^1$. Then it is clear
that the equation
$\gamma^2(z)\ = \ 1\ +\ \alpha(z)F(z)$
can be solved {\it in a neighborhood} $U$ of $z_0$. Thus we obtain that
$\beta(z)\ =\ \alpha(z)\gamma(z)^2\quad\text{for}\quad z\in U$
and so {\it the forms in $\CC$ represented by the functions $\alpha$ and 
$\beta$ become isometric, if we ignore the points of the circle outside $U$}.

This shows that despite the drastic global 
changes which can be performed on the 
function $\alpha$ and which do not influence the corresponding Hermitian 
form $\phi_\alpha$ in $\tca$, the behavior of the function $\alpha(z)$
{\it near the zeros} $\alpha(z)=0$ remains essentially the same. We will see in the next subsection, 
that a similar (and more precise!) statement can be made in a completely general situation.

\subheading{4.6. Excision} We will describe now a procedure of excision
for Hermitian forms in $\tca$, which is similar to the procedure of
excision on objects of the extended category, cf. \cite{Fa3}, 1.4. 

Suppose that the given Hermitian form $\alpha:A\to D(A)$ can be 
represented as an orthogonal sum 
$A\ =\ P\oplus Q, \quad P\perp Q,$ where $P$ and $Q$ are orthogonal with
respect to the form $\alpha$,
and the restriction $\alpha|_P$ is non-degenerate, i.e. $\alpha|_P: P\to D(P)$
is an isomorphism. Then the Hermitian form 
$$\beta=\alpha|_Q: Q\to D(Q)$$
produces discriminant form $\phi_\beta$ on $(\beta:Q\to D(Q)),$
which is obviously isometric to $\phi_\alpha$ in $\tca$.
We will say that {\it the form $\beta$ is obtained from $\alpha$ by excision
with respect to $P\subset A$.}

As an example of excision, which we will use frequently, 
consider the following construction. Suppose that we have a $\CC$-admissible
scalar product on $A$ and a self-adjoint operator $\alpha\in\hom_{\CC}(A,A)$.
It defines a Hermitian form in $\tca$, cf. 4.4. Let 
$$\alpha\ =\ \int_{-\infty}^\infty\lambda dE_\lambda\tag4-11$$
be the spectral representation of $\alpha$. Given $\epsilon >0,$ consider 
the decomposition $A\ =\ P\oplus Q$, where 
$$P\ =\ E_{-\epsilon}A\ + \ (1-E_\epsilon)A,\quad Q\ 
=\ (E_\epsilon  - E_{-\epsilon})A.\tag4-12$$ 
Then $\alpha$ maps $P$ into $P$ and induces an isomorphism $P\to P$. Also,
$P$ and $Q$ are orthogonal with respect to $\alpha$. Thus
we can make an excision with respect to $P$ to get a form $\alpha|_Q:Q\to Q$.

Note that the norm of $\alpha|_Q$ is less or equal than $\epsilon$.

The following two theorems will be our main technical results allowing to reduce
the questions about non-degenerate torsion forms in $\tca$ to questions
about (degenerate!) forms in the initial von Neumann category $ca$.

\proclaim{4.7. Theorem} Any non-degenerate Hermitian form in $\tca$
is isometric to a discriminant form of a Hermitian form in $\CC$.
\endproclaim

\proclaim{4.8. Theorem} Suppose that two Hermitian forms $\alpha:A\to D(A)$ 
and $\beta:B\to D(B)$ in category $\CC$ (with injective $\alpha$ and 
$\beta$) produce isometric discriminant forms $\phi_\alpha$ and $\phi_\beta$ 
in $\tca$. Then the forms $\alpha$ and $\beta$ admit excisions, which are
isometric as forms in $\CC$.
\endproclaim

\subheading{4.9. Proof of Theorem 4.7} Without loss of generality, we may assume that we are given a 
non-degenerate Hermitian form $\phi:\X\to \e(\X)$ in $\tca$ represented by
a commutative diagram
$$
\CD
(\alpha: A @>>> A)\\ 
@V{f^\ast}VV @VVfV\\
(\alpha: A @>>> A)
\endCD \tag4-13
$$
(a symmetric presentation) where we suppose that a $\CC$-admissible scalar
product on $A$ has been chosen and $\alpha$ is injective, self-adjoint
$\alpha^\ast =\alpha$, and has dense image. 

Suppose that $f$ in diagram (4-13) is an isomorphism. Then we have a 
discriminant form 
$$
\CD
(f\alpha: A @>>> A)\\ 
@V{\id}VV @VV{\id}V\\
(\alpha f^\ast: A @>>> A),
\endCD \tag4-14
$$
which is congruent to (4-13) because of the diagram
$$
\CD
(\alpha: A @>>> A)\\ 
@V{\id}VV @VVfV\\
(f\alpha: A @>>> A)\\
@V{\id}VV @VV{\id}V\\
(\alpha f^\ast: A @>>> A)\\  
@V{f^\ast}VV @VV{\id}V\\
(\alpha: A @>>> A).
\endCD\tag4-15
$$

Thus, to prove the Theorem, we have to show that, given diagram (4-13), 
representing an isomorphism in $\tca$, we can "replace" it by a diagram with
$f$ isomorphic. The "replacement" will mean certain "excision". 

As {\bf the first step}, we want to show that we may make $f$ injective. 
Let $B\subset A$ denote $\overline{\im (f^\ast)}$ and $B'\subset A$ denote 
$\overline{\im (f)}$; we consider
$B$ and $B'$ with the induced from $A$ scalar products. Denote by $i:B\to A$ 
and $i':B'\to A$ the inclusions, and $\pi:A\to B$ and $\pi':A\to B'$ will 
denote the orthogonal projections. 

We will also need the following map $j': B'\to A$. To define it, note, 
that the orthogonal complement to $B'$ coincides with $\ker(f^\ast)$ and 
the orthogonal complement to $B$ is $\ker(f)$. Since diagram (4-13) 
represents a monomorphism, using Proposition 1.6 of \cite{Fa3}, we find 
that $\alpha$ maps $\ker(f^\ast)$ {\it onto} $\ker(f)$. Therefore, given 
$x\in B'$ we may write $\alpha(x) = a + b$, where $a\in B$ and $b\in \ker(f)$. 
Now, as explained above, $b=\alpha(b')$ for some $b'\in\ker(f^\ast)$ and 
we set
$$j'(x) = x - b'.\tag4-16$$
The map $j'$ is clearly a monomorphism. 

Consider the torsion Hilbertian module
$$\Y = (\beta: B'\to B),\quad\text{where}\quad 
\beta = \pi\circ \alpha\circ i'.$$
We have the following morphism in $\tca$
$$
\CD
(\beta: B'@>>> B)\\
@V{j'}VV @VV{i}V\\
(\alpha: A @>>> A).
\endCD\tag4-17
$$
This implies that $\beta$ is injective and has dense image.
Using Proposition 1.6 of \cite{Fa3} we check that the kernel and the cokernel
of the morphism of $\eca$, represented by
(4-17) are trivial and so (4-17) represents an isomorphism in $\eca$. 

The initially given form $\phi$ on $\X$ induces a form $\psi: \Y\to \e(\Y)$ via 
isomorphism (4-17). We find that $\psi$ is given by the diagram
$$
\CD
(\beta: B'@>>> B)\\  
@V{g^\ast}VV @VV{g}V\\  
(\beta^\ast: B @>>> B'),
\endCD\tag4-18
$$
where $g={j'}^\ast\circ f\circ i:B\to B'$. We only have to check now 
that $g$ is
{\it injective}. But $g$ clearly coincides with the action of $f$
as the map $\overline{\im(f^\ast)} \to \overline{\im(f)}$. The last map 
is obviously injective on $\im(f^\ast)$; now the result follows because of 
finiteness assumptions, using  \cite{Fa3}, Proposition 2.4. Thus, we
may replace the initial form (4-13) by (4-18) where $g$ is injective.

As {\bf the second step}, we will assume that $f$ in (4-13) is injective 
and will make it isomorphism by performing further excisions.

Since diagram (4-13) represents an isomorphism, using Proposition 1.6 of
\cite{Fa3}, we find that the sequence
$$0\to A @>{\alpha\oplus - f^\ast}>>A\oplus A @>{(f,\alpha)}>>
A\to 0\tag4-19$$
is exact in $\eca$. Therefore it splits, since it consists of projective
objects. Thus, there exist morphisms
$\sigma:A\to A$ and $\delta: A \to A$, such that
$$\id_{A}\ =\ \sigma\alpha\ +\ \delta f^\ast.\tag4-20$$
This is the crucial point of the proof; intuitively, (4-20) means that 
{\it $f^\ast$ and $\alpha$ cannot be small at the
same places}. Choose $\epsilon > 0$ such that $\epsilon\cdot ||\sigma||<1$.
Using decomposition (4-12), we may construct a splitting
$$A =\ P \oplus\ Q,\quad P\perp Q,\tag4-21$$ 
such that $\alpha(P)=P$, $\alpha(Q)\subset Q$ and the norm
of the restriction of $\alpha$ on $Q$ is less than $\epsilon$, i.e.
$||\alpha|_{Q}||< \epsilon$.
Denoting by $i_{Q}: Q \to A$ and by $\pi_{Q}: A \to Q$ the corresponding 
inclusion and projection, we will have (using (4-20)) that the
morphism $\pi_Q \delta f^\ast i_{Q}:Q\to Q$ is an isomorphism. Therefore,
$f^\ast i_Q:Q\to A$ is injective with closed image.

Let's show that the (injective) map $fi_Q:Q\to A$ also has closed image.
Since (4-13) represents an epimorphism in $\eca$, we know that 
$\im(\alpha) + \im(f) = A$ (cf. Proposition 1.6 of \cite{Fa3}). We have
$A=\overline{f(P)} \oplus \overline{f(Q)}$ and thus we obtain 
(since we have already proved that $f^\ast(Q)\subset A$ is closed) that
$\alpha(f^\ast(Q)) + f(Q)=\overline{f(Q)}$. 
Thus, since $\alpha(f^\ast(Q))\subset f(Q)$, we get $f(Q)=\overline{f(Q)}$.

We now want to show that the map $\pi_Q : Q'\to Q$ is an isomorphism; it 
is clearly equivalent to the statement that $\pi_Q f^\ast i_Q: Q\to Q$ is an 
isomorphism. We have
$$\pi_Q f i_P = \alpha \pi_Q f^\ast \alpha^{-1} i_P$$
(because $f\alpha = \alpha f^\ast$ and $\alpha$ preserves $P$ and $Q$ and
is an isomorphism $\alpha: P\to P$). Therefore,
$$||\pi_Q f i_P|| \le C\epsilon$$
for some constant $C >0$ (since $||\alpha|_Q||<\epsilon$), and thus
$$||\pi_P f^\ast (q)|| \le C\epsilon\cdot ||q||\quad\text{for}\quad q\in Q.$$
On the other hand, since we know that $f^\ast :Q\to A$ is injective with 
closed image, holds
$$||f^\ast(q)|| \ge C'\cdot ||q||,\quad\text{for}\quad q\in Q$$
for some constant $C'>0$. Hence we obtain
$$
\aligned
||\pi_Q f^\ast (q)|| & \ge ||f^\ast(q)|| - ||\pi_P f^\ast(q)|| \ge \\
& \ge (C' -C\epsilon)\cdot ||q||,\quad\text{for}\quad q\in Q. 
\endaligned \tag4-22  
$$
Note that the constants $C$ and $C'$ may be assumed independent of 
$\epsilon$ and of the decomposition $A= P+Q$ and so for small enough 
$\epsilon>0$ we will have $\pi_q f^\ast :Q \to Q$ isomorphism because of 
(4-22). 

Now we may finish the proof as follows. Denote $Q' = F(Q)$. The diagram
$$
\CD
Q @>{f\alpha i_Q}>> Q'\\
@V{i_Q}VV @VV{f^{-1}}V\\
A @>{\alpha}>> A 
\endCD\tag4-23
$$
represents an isomorphism in $\eca$ (apply Proposition 1.6 of \cite{Fa3}), 
and the induced Hermitian form on $(f\alpha i_Q: Q \to Q')$ is obviously
given by the diagram
$$
\CD
Q @>{f\alpha i_Q}>> Q'\\
@V{\pi_{Q'}i_Q}VV @VV{\pi_Q i_{Q'}}V\\
Q' @>{(f\alpha i_Q)^\ast}>> Q 
\endCD\tag4-24
$$
We have shown above that $\pi_Q i_{Q'}: Q'\to Q$ is an isomorphism. The 
result now follows from the argument in the very beginning of the proof.

\subheading{4.10. Proof of Theorem 4.8} 
Suppose that some $\CC$-admissible scalar products on $A$ and 
on $B$
are chosen so that we may consider $\alpha:A\to A$  and $\beta:B\to B$ as
self-adjoint operators, injective, with dense images. An isometry between 
the discriminant forms $\phi_\alpha$ and $\phi_\beta$ is then represented 
by commutative diagram of the form
$$
\CD
A @>\alpha>>A\\
@VgVV   @VVfV\\
B @>>\beta >B.
\endCD\tag4-25
$$
The fact, that it is an isometry means that it is an isomorphism in $\eca$
and also the composition morphism given by the diagram
$$
\CD
A @>\alpha>>A\\
@VgVV   @VVfV\\
B @>>\beta >B\\
@V{f^\ast}VV   @VV{g^\ast}V\\
A @>{\alpha}>>A,
\endCD\tag4-26  
$$
equals $\phi_\alpha$ (defined as (4-9)). Therefore, there exists
a morphism $F: A\to A$ of $\CC$ such that 
$$g^\ast\circ f = \id_A + \alpha\circ F,\quad 
f^\ast\circ g = \id_A + F\circ \alpha.\tag4-27$$

Now, choose $\epsilon>0$ such that $\epsilon \cdot ||F|| < 1$, and perform
excision (4-12) on $(\alpha: A\to A)$ with respect to this $\epsilon$. 
As the result we will obtain $(\alpha|_Q: Q\to Q)$ and (4-25) produces
the isometry
$$
\CD
Q @>{\alpha|_Q}>> Q\\
@V{g|_Q}VV @VV{f|_Q}V\\
B @>>\beta>B.
\endCD\tag4-28 
$$
From the second equation in (4-27) we see that $g|_Q$ is injective 
with closed image (this is similar to the argument used in the proof of 
Theorem 4.7).

To simplify our notations, we will now assume that the initial isometry
(4-25) has the property
$$||\alpha||\cdot ||F|| < 1 - \delta\quad\text{for some small}
\quad \delta>0,\tag4-29$$
where $F$ satisfies (4-27). Then, as we have already mentioned, $g$ is
injective with closed image. We will denote $B'=g(B)\subset A$.

We will finish the proof of the Theorem by showing:
\roster
\item"{(i)}" {\it The Hermitian form $\alpha$ is congruent to the form
$\beta'$ induced by $\beta$ on $B' \subset B$;}
\item"{(ii)}" {\it The form induced by $\beta$ on ${B'}^\perp$, the 
orthogonal complement to $B'$ with respect to $\beta$, is non-degenerate.}
\endroster

Indeed, (i) would imply that $\beta'$ is {\it weakly non-degenerate}, i.e.
the morphism $\beta': B'\to D(B')$ is injective; therefore, 
$A = B'\oplus {B'}^\perp$, and so using (ii) we obtain that $\alpha$ {\it is 
isometric to an excision of} $\beta$.

On the other hand, once (i) is obtained, (ii) would follow from the
fact that (4-25) is an isomorphism in $\eca$. So, we only need to prove
(i).

Consider the form $A$ induced by $g$ from $\beta$, i.e. $g^\ast \beta g$.
From (4-26) and (4-27) we obtain
$$g^\ast \beta g = \alpha f^\ast g = \alpha + \alpha F \alpha.\tag4-30$$

Our aim (i) will be achieved if we will show that the right hand side of 
(4-30) can be represented
in the form 
$$\alpha\ +\ \alpha F\alpha \ = \ h^\ast \alpha h,
\quad\text{where} \quad h:A\to A\tag4-31$$
is an $\CC$-isomorphism.

Let us represent $\alpha$ in the form
$\alpha\ =\ s\gamma^2,$
where the morphisms $s,\gamma:A\to A$ satisfy
$$s^2\ =\ 1, \quad s^\ast=s,\quad \gamma^\ast=\gamma,\quad \gamma>0,\quad
\gamma s=s\gamma.$$
The morphism $\gamma$ is the positive square root of $|\alpha|$, and $s$ is
the {\it sign} of $\alpha$. The operators $\gamma$ and $s$ are obtained by
$$s = \int_{-\infty}^{\infty} \sign(\lambda)dE_\lambda,\quad
\gamma = \int_{-\infty}^{\infty} |\lambda|dE_\lambda,$$
where $E_\lambda$ are the spectral projectors of the spectral decomposition 
of $\alpha$.
Note, that all the spectral projectors and the operators $s$ and $\gamma$,
given by the above integrals,
are morphisms of $\CC$ since $\hom_{\CC}(A,A)$ is a von Neumann algebra, 
and so we may use the general results of \cite{Di}.

Then we have 
$\alpha\ +\ \alpha F\alpha \ 
=\  s\gamma(1+\gamma Fs\gamma)\gamma.$

Now, let us find $\CC$-isomorphisms
$h, h_1, h_2 :\ A\ \to A,$
such that 
$$
\aligned
1+\gamma Fs\gamma \ =&\ h_1^2,\\
1+\gamma^2Fs\ =&\ h_2^2,\\
1+Fs\gamma^2\ =&\ h^2.
\endaligned\tag4-32
$$
We will show that we may construct $h$, $h_1$ and $h_2$ such that the 
following identities will hold
$$s h_2\ =\ h^\ast s,\tag4-33$$
$$\gamma h_1\ =\ h_2 \gamma,\tag4-34$$
$$\gamma h \ =\ h_1 \gamma.\tag4-35$$
To construct $\delta$ we observe that the spectrum of $1+\gamma Fs\gamma$
is contained (because of (4-29)) inside the circle
$\Gamma_\delta\ =\ \{\lambda; |\lambda-1|=1-\delta\}.$
Thus we may use the functional calculus to define 
$$h_1 \ =\ \frac{1}{2\pi i}\int_{\Gamma_\delta} 
\frac{{\lambda}^{1/2}}{\lambda-(1+\gamma F s \gamma)}d\lambda,\tag4-36$$
where $\lambda^{1/2}$ denotes the branch of the square root which is obtained 
by cutting along the negative real axis; we will use the fact that it 
commutes with the complex conjugation. We define the operators $h_2$ and 
$h$ by the similar integrals
$$h_2 \ =\ \frac{1}{2\pi i}\int_{\Gamma_\delta} 
\frac{{\lambda}^{1/2}}{\lambda-(1+\gamma^2 F s)}d\lambda\tag4-37$$
and
$$h \ =\ \frac{1}{2\pi i}\int_{\Gamma_\delta} 
\frac{{\lambda}^{1/2}}{\lambda-(1+ F s \gamma^2)}d\lambda.\tag4-38$$
One easily checks that the 
formulae (4-33), (4-34), (4-35) are clearly satisfied. For example to
prove (4-34) one observes that
$$\gamma\cdot(\lambda-(1+\gamma F s \gamma))^{-1} = 
(\lambda-(1+\gamma^2 F s))^{-1}\gamma.$$ 

Now we have using (4-33), (4-34), (4-35)
$$ 
\aligned
\alpha\ +\ \alpha F\alpha \ =&\ s\gamma h_1^2\gamma\\  
=& \ s(\gamma h_1)(h_1 \gamma)\ \\
=& \ s(h_2 \gamma)(\gamma h)\ \\
=& \ h^\ast s\gamma^2 h\  \\
=& \ h^\ast\alpha h.
\endaligned\tag4-39
$$
The fact that $h$ is an isomorphism follows by comparing the 
equation $1+F\alpha = h^2$ (cf. (4-33)) with inequality (4-29).

This completes the proof. \qed

\heading{\bf \S 5. Positively and negatively definite torsion Hermitian forms}\endheading

\subheading{5.1. Definition}
{\it A non-degene\-rate Hermitian form $\phi: \X\to \e(\X)$ in $\tca$ will be called 
positively (negatively) definite if it can be represented as the 
discriminant form $\phi_\alpha$ of a positively (negatively) defined 
Hermitian form} $\alpha:A\to D(A)$ in $\CC$.

Equivalently, a torsion Hermitian form 
$\phi: \X\to \e(\X)$ is positively definite, if and only if any Hermitian
form $\alpha: A\to D(A)$ in $\CC$ which produces $\phi$ as the discriminant
form $\phi = \phi_\alpha$, 
can be represented as an orthogonal sum $\alpha = \alpha_1 \perp
\alpha_2$, where $\alpha_1$ is non-degenerate and $\alpha_2$ is positively
definite; this follows from Theorem 4.8.

\proclaim{5.2. Lemma} Any non-degenerate torsion Hermitian form 
$\phi:\X\to \e(\X)$ in $\tca$, which at the same time is positively and negatively 
definite, is trivial: $\X=0$.\endproclaim
\demo{Proof} Suppose that a non-degenerate torsion Hermitian form   
can be represented as the discriminant form of a positively definite form 
$\alpha:A\to D(A)$ and of a negatively definite form $\beta:B\to D(B)$.
Then by Theorem 4.8, $\alpha$ and $\beta$ have isometric excisions.
Since any excision of a positively (negatively) defined form is positively
(negatively, correspondingly) defined, and since a positively definite form
can be isometric to a negatively definite form if and only if the both forms 
are trivial, we obtain that the forms $\alpha$ and $\beta$ admit trivial
excisions. Hence the maps $\alpha:A\to D(A)$ and $\beta:B\to D(B)$ are
isomorphisms and so $\X=0$. $\square$\enddemo

\proclaim{5.3. Theorem} Any non-degenerate torsion Hermitian form 
$\phi:\X\to \e(\X)$ in $\tca$ can be represented as an orthogonal sum
$$\phi\ =\ \phi_+\ \perp \ \phi_-,\tag4-41$$
where $\phi_+:\X_+\to \e(\X_+)$ is a positively definite and $\phi_-:\X_-\to \e(\X_-)$ 
is a negatively definite torsion Hermitian forms.
\endproclaim

\demo{Proof} It follows directly from Theorem 4.7.
Namely, given a non-degenerate torsion Hermitian form 
$\phi:\X\to D(\X)$ in $\tca$, by Theorem 4.7 it can be represented
as the discriminant form $\phi_\alpha$ of a Hermitian form
$\alpha\in\hom_{\CC}(A, D(A))$. If we chose a $\CC$-admissible scalar
product on $A$, then $\alpha$ can be identified with a self-adjoint operator
$\alpha: A\to A$. Using the spectral theorem, we can represent $\alpha$ as 
$\alpha_+\oplus \alpha_-$ where $\alpha_\pm:A_\pm\to A_\pm$ is a positive 
(negative) operator $\alpha_\pm\in\hom_{\CC}(A,A)$. Then $\phi\ =\ 
\phi_+\oplus \phi_-$, where $\phi_\pm$ is the discriminant form of 
$\alpha_\pm$. \qed
\enddemo

We will prove later uniqueness of the decomposition of Theorem 5.3 (cf. Theorems 7.7 and
7.4), but we will impose some additionally assumption on the category (superfiniteness).

\heading{\bf \S 6. Metabolic and hyperbolic forms}\endheading

In this section we will assume that $\CC$ is a finite von Neumann category, and thus $\tca$ is an
abelian category, cf. \cite{Fa3}.

Recall a general terminology concerning Hermitian forms. 
A form $\phi: \X \to \e(\X)$ in $\tca$
is called {\it metabolic} if there is an inclusion $i: \Y\to\X$ (in the sense of the category theory)
such that $\Y$ coincides
with its annihilator $\Y^\perp$, i.e. with the kernel of
$$\X @>\phi>> \e(\X) @>{\e(i)}>> \e(\Y).\tag6-1$$
Such subobject $\Y\subset \X$ with $\Y = \Y^\perp$ is called {\it a metabolizer}.

A form $\phi: \X \to \e(\X)$ is called {\it hyperbolic} if it has a metabolizer which is
a direct summand.

We will see that torsion Hermitian forms in $\tca$ are all metabolic, and, moreover, in
many cases the metabolizer is unique. We will also see that they are rarely hyperbolic. 

\proclaim{6.1. Proposition} Any non-degenerate torsion Hermitian form $\phi: \X \to \e(\X)$ in $\tca$
 is metabolic.
\endproclaim
\demo{Proof} Because of Theorem 5.3, it is enough to prove metabolicity of positively
defined torsion forms. By Theorem 4.7 and observations made in 4.4, we may assume that
the given form $\phi$ is represented by diagram (4-9), where $\alpha: A\to A$ is a
self-adjoint positive operator, $\alpha\in \hom_{\CC}(A,A)$. We may
find a positive square root of $\alpha$, i.e. $\beta \in \hom_{\CC}(A,A)$
with $\beta^2 =\alpha$, $\beta^\ast = \beta$, $\beta > 0$. Denote $\Y = (\beta: A\to A)$,
and let $i: \Y \to \X$ be given by the diagram
$$
\CD
(A @>{\beta}>> A)\\
@VV1V      @V{\beta}VV\\
(A @>\alpha>> A).
\endCD
$$
Then the composition (6-1) equals
$$
\CD
(A @>{\alpha}>> A)\\
@VV{\beta}V      @V{1}VV\\
(A @>{\beta}>> A)
\endCD
$$
and using Proposition 1.6 of \cite{Fa}, one easily checks that its kernel coincides
with $i: \Y \to \X$. \qed
\enddemo

\proclaim{6.2. Proposition} Let $\phi: \X \to \e(\X)$ be a non-degenerate Hermitian form, which
is positively (or negatively) definite. Then in $\X$ there is a unique metabolizer.
\endproclaim
\demo{Proof} Suppose that $\phi: \X \to \e(\X)$ is presented as the discriminant form
of a positively definite Hermitian form $\alpha: A\to A$.
As in the proof of 6.1 we suppose that $A$ has a specified admissible scalar product and
$\alpha \in \hom_{\CC}(A,A)$ is self-adjoint and positive.
Then any inclusion $i: \Y \to \X$ in $\tca$
can be represented by a diagram
$$
\CD
(A @>{\beta}>> A)\\
@V1VV @VV{\gamma}V\\
(A @>{\alpha}>> A),
\endCD
$$
where $\alpha = \gamma \beta$. Using Proposition 1.6 of \cite{Fa}
we observe that $\Y \subset \Y^\perp$ if and only if
$\gamma$ can be represented as $\gamma = \beta^\ast \delta$ for some
$\delta \in \hom_{\CC}(A,A)$.
Thus we have
$\alpha = \beta^\ast \delta \beta.$
Now, $\alpha^\ast =\alpha$ implies $\delta^\ast = \delta$. Also, $\delta$ is clearly
positive. 
The inclusion $\Y^\perp \to \X$ is represented by
$$
\CD
(A @>{\gamma}>> A)\\
@V1VV @VV{\beta}V\\
(A @>{\alpha}>> A)
\endCD
$$
and the inclusion $\Y \to \Y^\perp$ is represented by
$$
\CD
(A @>{\beta}>> A)\\
@V1VV @VV{\delta}V\\
(A @>{\gamma}>> A).
\endCD
$$
Again, using Proposition 1.6 of \cite{Fa} we find that $\Y = \Y^\perp$ (i.e. $\Y$ is a
metabolizer) if and only if
$\delta$ is an isomorphism in $\CC$. Now we may write $\delta = \delta_1^2$ with
$\delta_1^\ast = \delta_1$
and $\delta_1 >0$. Therefore we see that our metabolizer $\Y$ coincides with the one
constructed in Proposition 6.1 (which corresponds to $\delta_1\beta$ taken instead of
$\beta$). \qed

\enddemo

\heading{\bf \S 7. Superfinite von Neumann categories}\endheading

In this section we introduce an additional finiteness assumption on the initial von Neumann
category $\CC$, which we call {\it superfiniteness}. We classify completely torsion Hermitian
forms in $\tca$, assuming that $\CC$ is superfinite. More precisely, we prove 
(cf. Theorems 7.4 and 7.7)
that any such
form determines canonically two torsion modules $\X_+$ and $\X_-$
(the positive and the negative parts)
and, conversely, the congruence class of the form is determined by the isomorphisms types of
$\X_+$ and $\X_-$. 

\proclaim{7.1. Proposition} For a finite von Neumann category $\CC$ the following
properties are equivalent:
\roster
\item"{(i)}"  any epimorphism $\X \to \X$, where $\X\in\ob(\tca)$, is an isomorphism;
\item"{(ii)}" any monomorphism $\X \to \X$, where $\X\in\ob(\tca)$, is an isomorphism;
\item"{(iii)}" in any commutative diagram of the form 
$$
\CD
A @>{\alpha}>> A\\
@VgVV @VVfV\\
A @>{\alpha}>> A,
\endCD\tag7-1
$$ 
in $\CC$, where $\alpha$ is injective and $f$ is an isomorphism in $\CC$, the morphism $g$ is also
an isomorphism in $\CC$.
\endroster
\endproclaim
\demo{Proof} Equivalence of (i) and (ii) easily follows from existence of the duality
$\e$ in $\tca$, cf. 3.3.

Suppose that (i) and (ii) are satisfied. Then any diagram (7-1) can be viewed as a morphism
$\phi : \X \to \X$ of the torsion subcategory $\tca$ where $\X = (\alpha: A \to A)$. If $f$
is an isomorphism then $\phi$ is an epimorphism, and thus it is an isomorphism by (i).
This implies that $g$ is an isomorphism in $\CC$ (as easy follows from Proposition 1.6 of
\cite{Fa3}).

Conversely, suppose that we know that (iii) is always satisfied. Using the arguments used in the
proof of Theorem 4.7, we obtain that any epimorphism $\X \to \X$ can be represented by a
diagram of form (7-1), where $\X = (\alpha: A \to A)$, $\alpha$ is injective and $f$ is
an isomorphism. Then (iii) implies that $g$ is an isomorphism as well, and thus the
given epimorphism $\X \to \X$ is an isomorphism. Therefore, (iii) implies (i).
\qed
\enddemo

\subheading{7.2. Definition} A finite von Neumann category $\CC$ will be called
{\it superfinite} if the equivalent conditions of Proposition 7.1 hold for $\CC$.

The diagram (7-1) can be viewed as follows. We may think of $\alpha$ as defining a new
scalar product on $A$ and then the operators $f$ and $g$ describe "the same" operator,
which is bounded with respect to two different norms. Such operators are studied in \cite{GK},
for example, where they called bibounded. I am grateful to V. Matsaev, who explained this to
me. Our condition of superfiniteness (7-1) is equivalent to the statement that the
spectra of $f$ and $g$ coincide.

As example in \cite{GK}, chapter 5, \S 6 show, the
category of Hilbert spaces and bounded linear maps is not superfinite.

\subheading{7.3. Example} Consider the von Neumann category $\CC$ described as Example 2
in section 2.6. A diagram of type (7-1) in this category looks as follows.

We are given a locally compact Hausdorff space $Z$ with a positive Radon measure $\mu$,
and the Hilbert space $A$ is the direct integral (2-5) of a measurable
field of finite dimensional Hilbert spaces $\xi \to \H(\xi)$, where $\xi \in Z$, cf.
\cite{Di}, part II.
We will assume that $\H(\xi)$ is nonzero for any $\xi \in Z$.
The morphisms $\alpha$, $f$, and $g$ of (7-1) are given by decomposable linear maps
(2-6) constructed out of essentially bounded measurable fields of linear maps
$\H(\xi) \to \H(\xi)$, cf. 2.7. We will denote the corresponding measurable fields of
linear maps by
$$T_\alpha(\xi) ,\quad T_f(\xi),\quad T_g(\xi) : \ \H(\xi) \to \H(\xi), \qquad \xi \in Z.\tag7-2$$
Since $\alpha$ is supposed to be injective, we know that $T_\alpha(\xi)$ is an isomorphism
for almost all $\xi\in Z$. Since we have the identity
$$ T_g(\xi) = T_\alpha (\xi)^{-1}\circ T_f (\xi) \circ T_\alpha (\xi)\tag7-3$$
for almost all $\xi$, we obtain that $T_f (\xi)$ and $T_g (\xi)$ have the same spectrum
for almost all $\xi$. We assumed that $f$ is an isomorphism, so
the field $\xi \mapsto T_f (\xi)^{-1}$ is essentially bounded and therefore the
spectrum of the operator $T_g (\xi)$ does not approach zero. More precisely, we obtain
that there exists $r > 0$ such
that the spectrum of $T_g (\xi)$ does not intersect the circle $\{\lam; |\lam | \le r\}$
for almost
all $\xi\in Z$. Also, the norm of $T_g (\xi)$ has a uniform upper bound (since $g$ is
bounded). Thus we obtain that
$$||T_g (\xi)^{-1}|| \le ||T_g (\xi)||^{-1} \le 1/r$$
for almost all $\xi \in Z$ and therefore 
$g$ is an isomorphism in category $\CC$.

This shows that the von Neumann category of this example is superfinite.

One may conjecture that {\it any finite von Neumann category is superfinite}.

Here is the first result, whose proof use superfiniteness.

\proclaim{7.4. Theorem} If $\CC$ is a superfinite von Neumann category then any two
non-degene\-rate positively definite torsion Hermitian forms
$\phi: \X_1 \to \e(\X_1)$ and $\psi: \X_2 \to \e(\X_2)$ are congruent if and only if the
underlying
objects $\X_1$ and $\X_2$ are isomorphic as objects of $\tca$.\endproclaim

\demo{Proof} It is clear that congruence of $\phi$ and $\psi$ implies isomorphism
between $\X_1$ and $\X_2$ and 
we only need to prove the inverse statement.

Without loss of generality, using Theorems 4.7 and 4.8,
we may assume that the forms $\phi$ and $\psi$ are given as follows.
Let $A$ be an object of $\CC$ with a particular choice of a $\CC$-admissible scalar
product $\<\ ,\ \>$ on $A$. 
Let $\alpha_\phi, \alpha_\psi: A \to A$ be the self-adjoint positive operators,
$\alpha_\phi, \alpha_\psi\in \hom_{\CC}(A,A)$ such that the forms $\phi$ and $\psi$ are the
discriminant forms corresponding to $\alpha_\phi$ and $\alpha_\psi$ correspondingly as explained in
subsection 4.4. Then we have $\X_1 = (\alpha_\phi: A \to A)$ and $\X_2 = (\alpha_\psi: A \to A)$,
and since we know that $\X_1$ and $\X_2$ are isomorphic, we obtain (using Proposition 3.6
of \cite{Fa3})
that there is a commutative diagram
$$
\CD
A @>{\alpha_\phi}>> A\\
@VfVV    @VVgV\\
A @>{\alpha_\psi}>> A,
\endCD\tag7-4
$$
where $f$ and $g$ are isomorphisms in $\CC$. 
Now we have $g\alpha_\phi = \alpha_\psi f$ and therefore $f^\ast \alpha_\psi = \alpha_\phi g^\ast$
and $f^\ast g \alpha_\phi = \alpha_\phi g^\ast f$. To be able to apply Lemma 7.5 below to
conclude that the spectrum of $f^\ast g$ is contained in an interval of the
positive real axis of the form $(\epsilon, N)$ for some small $\epsilon >0$
and for some large $N>0$. We only have to note that
$\alpha_\phi \ge 0$ and $(f^\ast g)\alpha_\phi = f^\ast \alpha_\psi f \ge 0$ using the
fact that $\phi$ and $\psi$ are positively definite.

Define
$$k = \frac{1}{2\pi i}\int_\Gamma \frac{\sqrt \lambda}{\lambda - f^\ast g}
d\lambda,\tag7-5$$
where $\Gamma$ is the boundary of the rectangle with vertices
$\epsilon -i\delta$, $N-i\delta$, $N+i\delta$ and $\epsilon + i\delta$
for some $\delta >0$. Here we use the branch of the $\sqrt{\lambda}$, 
obtained by cutting along the negative real axis, {\it which commutes with the
complex conjugation}. Then, using the 
functional calculus \cite{DS}, chapter VII, \S 3, we have 
$$k^2 = f^\ast g\tag7-6$$
and, in particular, $k$ is a $\CC$-isomorphism; moreover, 
$$k \alpha_\phi = \alpha_\phi k^\ast.\tag7-7$$   
To prove the last formula observe that 
$$k^\ast = \frac{1}{2\pi i}\int_\Gamma \frac{\sqrt \lambda}
{\lambda - g^\ast f}d\lambda\tag7-8$$
and now (7-7) follows from the identity
$$(\lambda - f^\ast g)^{-1} \alpha_\phi = \alpha_\phi (\lambda - g^\ast f)^{-1}.$$
Thus we obtain
$$f^\ast \alpha_\psi f = f^\ast g \alpha_\phi = k^2 \alpha_\phi = k \alpha_\phi k^\ast,\tag7-9$$
which shows that $\phi$ and $\psi$ are congruent.

Therefore the torsion forms $\phi$ and $\psi$ are congruent.
This completes the proof. \qed
\enddemo

Here is the Lemma which we used in the above proof.

\proclaim{7.5. Lemma} Suppose that $\CC$ is a superfinite von Neumann category and
$\alpha$ and $\beta$ are operators
acting in an object $A \in \ob(\CC)$, where $\alpha, \beta \in \hom_{\CC}(A, A)$.
Suppose that we have specified an admissible scalar product $\<\ ,\ \>$ on $A$ such that
$\beta$ and $\beta\alpha$ are
self-adjoint and positive. Then the spectrum of $\alpha$ is contained in the non-negative
part of the real axis.
\endproclaim

\demo{Proof} I am very thankful to V. Matsaev, who explained to me
that this statement {\it is not true} in the category of Hilbert
spaces: he showed to me examples, constructed in \cite{GK}, chapter 5, \S 6,
which do not satisfy Lemma 7.5. Here is the only
place, where we will use the
assumption of superfiniteness of our category $\CC$.

Denote by $\<\ ,\ \>_\beta$ the scalar product on $A$ determined by
$\beta$, where $\<h,h'\>_\beta = \<\beta(h),h'\>$, (for $h, h'\in A$)
and by $||\cdot||_\beta$ 
the corresponding norm on $A$.
We have for $h\in A$, $h\ne 0$,
$$\<\alpha(h),h\>_\beta \, \, >, \, \,  0\tag7-10$$
i.e. $\alpha$ is self-adjoint and positive with respect to $\beta$-scalar product.
If $\beta$ is invertible, then the Lemma follows from the well-known fact
that the positive operator in Hilbert space has non-negative spectrum.

Now we use Lemma 7.6 below to conclude that
there exists a bounded linear operator $f: A \to A$ with $f\in \hom_{\CC}(A, A)$
such that the following digram
$$
\CD
A @>{\beta^{1/2}}>> A  @>{\beta^{1/2}}>> A\\
@V{\alpha}VV @VfVV @VV{\alpha^\ast}V\\
A @>{\beta^{1/2}}>> A  @>{\beta^{1/2}}>> A
\endCD\tag7-11
$$
is commutative. From this 
commutative diagram we see that $f$ is self-adjoint. Thus the spectrum of
$f$ belongs to the positive real axis. Therefore, for any $\lam\in \C$, which is not real and
positive, $f-\lam$ is invertible. Therefore, we may apply the definition of superfiniteness
(and Proposition 7.1) to the diagram
$$
\CD
A @>{\beta^{1/2}}>> A\\
@V{\alpha -\lam}VV @VV{f-\lam}V\\
A @>{\beta^{1/2}}>> A,
\endCD
$$
to conclude that $\alpha -\lam$ is invertible. Therefore, we obtain that
the spectrum of $\alpha$ is contained in the positive real axis. \qed
\enddemo

\proclaim{7.6. Lemma} Given a commutative diagram
$$
\CD
A @>{\alpha^2}>> A\\
@VfVV @VVgV\\
A @>{\alpha^2}>> A
\endCD
$$
in a finite von Neumann category, where $A$ is supplied with an admissible scalar product
and $\alpha$ is self-adjoint and positive. Then this diagram can be completed to a commutative
diagram
$$
\CD
A @>{\alpha}>> A  @>{\alpha}>> A\\
@VfVV @VhVV    @VVgV\\
A @>{\alpha}>> A  @>{\alpha}>> A,
\endCD
$$
where $h\in \hom_{\CC}(A,A)$.\endproclaim

Note that $h$ in the above diagram is clearly unique since $\alpha$ has dense image
(as follows from finiteness of $\CC$).

\demo{Proof} Existence of a bounded linear map $h$ such the above diagram is commutative,
is equivalent to Lemma 1.1 of \cite{DLS}. (Note that this Lemma of \cite{DLS} is a
generalization of a well-known theorem of M.G. Krein,
cf. also \cite{L}, Theorem I,  or \cite{D}).

In fact, in order to interpret the above mentioned results of \cite{DLS} in this form,
we note that
we may identify the inclusion of $A$ into its completion with respect to the $\alpha$-norm
$||x||_{\alpha} = ||\alpha(x)||$
with the morphism $\alpha: A \to A$ and after this identification the cited result
can be used immediately.

We claim now that the constructed operator $h$ must belong to the von Neumann algebra
$\hom_{\CC}(A, A)$. To show this, denote by
$\eta_\mu$, where $\mu > 0$, the following function of a real variable $\lam > 0$
$$
\eta_\mu(\lam) = \cases
1, \quad\text{for}\quad \lam < \mu\\
\lam,\quad\text{for}\quad \lam \ge \mu.
\endcases\tag7-12
$$
Note that the function $\lam \mapsto \eta_\mu(\lambda)^{-1}$ is bounded on the
spectrum of $\alpha$.
Now we see that $h$ is the weak limit
$$
h = \lim_{\mu\to 0} \ \alpha\circ f \circ
(\int_0^\infty \eta_\mu(\lam)^{-1} dE_\lam)\tag7-13$$
of the operators belonging to $\hom_{\CC}(A, A)$,
where $E_\lam$ is the spectral projection from the spectral decomposition of $\alpha$.
\qed
\enddemo

The following statement is the main result of this section. Together with Theorem 5.3
it gives a complete classification of torsion Hermitian forms.

\proclaim{7.7. Theorem} Let $\CC$ be a superfinite von Neumann category. Let
$\phi:\X \to \e(\X)$ and $\psi: \Y\to \e(\Y)$ be two torsion Hermitian forms in $\tca$.
Let $\X =  \X_+ \perp \X_- $ and $\Y =  \Y_+ \perp \Y_-$ be the decompositions into the 
positively and negatively definite parts, which exist
as asserted by Theorem 5.3. Then the initial forms $\phi$ and $\psi$ are congruent if and only if
$\X_+$ is isomorphic to $\Y_+$ and $\X_-$ is isomorphic to $\Y_-$ as objects of $\tca$.
\endproclaim

We will first prove the following Lemma:

\proclaim{7.8. Lemma} Suppose that $\CC$ is a superfinite von Neumann category. If
$\phi : \X \to \e(\X)$ is a non-degenerate positive torsion Hermitian form and
$\psi : \X \to \e(\X)$ is a positive torsion Hermitian form then
$\phi + \psi : \X \to \e(\X)$ is non-degenerate and positive.
\endproclaim
\demo{Proof} We defined in \S 4 the notion of positivity only for non-degenerate
Hermitian forms.
The form $\psi$ above in not supposed to be non-degenerate and its positivity means the
following. If $\CClK$ denotes the kernel of $\psi : \X \to \e(\X)$ then $\psi$ naturally
defines a non-degenerate form on $\X/\CClK$. We will say that $\psi$ is positive iff this
non-degenerate form on $\X/\CClK$ is positive.

Suppose that $\phi$ is given as the discriminant form of $\alpha: A\to A$, cf. Theorem
4.7. Here $A$ is an
object of $\CC$ with a choice of an admissible scalar product and $\alpha\in \hom_{\CC} (A, A)$
is self-adjoint and positive. The morphism $\phi$ is therefore represented by the diagram
$$
\CD
(A @>{\alpha}>> A)\\
@V1VV @VV1V\\
(A @>{\alpha}>> A).
\endCD
$$
The other form $\psi : \X \to \e(\X)$ is represented by a diagram
$$
\CD
(A @>{\alpha}>> A)\\
@V{\beta^\ast}VV @VV{\beta}V\\
(A @>{\alpha}>> A)
\endCD
$$
(the symmetric representation), cf. Lemma 4.2. Here $\alpha\beta^\ast = \beta\alpha$ and
so $\beta\alpha$ is self-adjoint. Now,
because of positivity of $\psi$ we may choose $\beta$ so that $\beta\alpha$ is non-negative.
The form $\phi + \psi$ is represented by the diagram
$$
\CD
(A @>{\alpha}>> A)\\
@V{1+\beta^\ast}VV @VV{1+\beta}V\\
(A @>{\alpha}>> A)
\endCD
$$
and the Lemma will be proven if we will show that $1+\beta$ is an isomorphism in $\CC$.

By Lemma 7.6 there exists $g\in \hom_{\CC}(A, A)$ such that the diagram
$$
\CD
A @>{\alpha^{1/2}}>> A @>{\alpha^{1/2}}>> A\\
@V{\beta^\ast}VV @VgVV @VV{\beta}V\\
A @>{\alpha^{1/2}}>> A @>{\alpha^{1/2}}>> A
\endCD
$$
commutes. We have $g \alpha^{1/2} = \alpha^{1/2} \beta^\ast$ and
$\beta \alpha^{1/2} = \alpha^{1/2} g$ which imply that $g$ is self-adjoint.
Since $\beta\alpha = \alpha^{1/2}g\alpha^{1/2} > 0$, we have $g > 0$. Therefore
$1 + g : A \to A$ is an isomorphism in $\CC$. Now considering the diagram
$$
\CD
A @>{\alpha^{1/2}}>> A @>{\alpha^{1/2}}>> A\\
@V{1+\beta^\ast}VV @V{1+g}VV @VV{1+\beta}V\\
A @>{\alpha^{1/2}}>> A @>{\alpha^{1/2}}>> A
\endCD
$$
and using superfiniteness (cf. Proposition 7.1) we obtain that $\beta^\ast +1$ and $\beta +1$
are both isomorphisms. \qed
\enddemo

\subheading{7.9. Proof of Theorem 7.7} We only have to prove the uniquenece 
of the decompositions into the positive and negative parts; the rest follows from 
Theorems 5.3 and 7.4.
Consider a non-degenerate torsion Hermitian form
$\phi: \X \to \e(\X)$ in $\tca$ and a decomposition $\X = \X_+ \perp \X_-$ given by Theorem
5.3. Suppose
that $i: \Y \to \X$ is a monomorphism in $\tca$ such that the induced form $\phi|_\Y$
on $\Y$
(i.e. the form $\e(i)\circ \phi\circ i : \Y \to \e(\Y)$) is non-degenerate
and positively definite.
The inclusion $i$ followed by the canonical projections $\pi_\pm : \X \to \X_\pm$,
determines the morphisms $f_\pm = \pi_\pm \circ i : \Y \to \X_\pm$. Denote by $\psi_\pm$ the
forms on $\X$ induced on $\Y$ by the morphisms $f_\pm$. We have
$$\phi|_\Y = \psi_+ + \psi_-.$$
We know that $\phi|_\Y$ is non-degenerate and positively definite. Also we know that
$-\psi_-$ is positively definite. Applying Lemma 7.8, we obtain that $\psi_+ =
\phi|_\Y + (-\psi_-)$ is non-degenerate.
Therefore the morphism $f_+: \Y \to \X_+$ is a monomorphism in $\tca$.

It follows that if $\X = \X'_+ \perp \X'_-$ is another representation
(4-46) then $\X'_+$ can be imbeded into $\X_+$ and similarly $\X_+$ can be imbeded
into $\X'_+$. Now, we obtain that $\X_+$ is isomorphic to $\X'_+$, because of superfiniteness of
$\CC$. \qed

\proclaim{7.10. Theorem} Let $\CC$ be a superfinite von Neumann category. Let
$\phi: \X \to \e(\X)$ in $\tca$ be a non-degenerate
torsion Hermitian form in $\tca$. Then $\phi$ is hyperbolic if and only if
the torsion objects $\X_+$ and $\X_-$ are isomorphic in $\tca$.
\endproclaim

\demo{Proof}Suppose first that we know that $\X_+ \simeq \X_-$ in the decomposition of the
form $\phi$ into an orthogonal sum of a positive and a negative forms. Then we may identify
$\X_+ = \X_- = \Y$ and using Theorem 7.4 we may assume that $\phi = \psi \perp (-\psi)$,
where $\psi : \Y \to \e(\Y)$ is a non-degenerate positive form on $\Y$. Then we have the
following metabolizer
$i: \Y \to \X = \Y \oplus \Y$, where $i = i_1 + i_2$, the sum of two inclusions.
This metabolizer is clearly a direct summand. Thus, $\phi$ is hyperbolic.

Now we want to prove the converse statement.
Suppose that $\phi: \X \to \e(\X)$ is Hermitian and non-degenerate and let
$\Y \subset \X$ be a metabolizer, $\Y = \Y^\perp$, in $\X$ which is a
direct summand, $\X = \Y \oplus \CClZ$. $\phi$ is a morphism $\Y \oplus \CClZ \to
\e(\Y) \oplus \e(\CClZ)$,
which is an isomorphism, and we know that the induced morphism $\Y \to \e(\Y)$ vanishes.
Therefore we obtain:
\roster
\item the induced morphism $\alpha: \Y \to \e(\CClZ)$ is a monomorphism (since its
kernel would be a part of the kernel of $\phi$);
\item the other induced morphism $\alpha^\dagger: \CClZ \to \e(\Y)$ (cf. notation introduced
in 1.3) is also a monomorphism (non-triviality of its kernel would contradict the condition
$\Y = \Y^\perp$).
\endroster
Since any torsion object is isomorphic to its dual (not canonically), we obtain that
$\Y$ can be imbeded into $\CClZ$ and conversely, $\CClZ$ can be imbeded into $\Y$. Using
superfiniteness of $\CC$, we obtain that $\Y$ and $\CClZ$ are isomorphic and, moreover,
the morphisms $\alpha: \Y \to \e(\CClZ)$ and $\alpha^\dagger: \CClZ \to \e(\Y)$
induced by $\phi$, are isomorphisms.

Denote by $\gamma: \CClZ \to \e(\CClZ)$ the morphism induced by $\phi$.
Note that, we may assume without loss of generality, that $\gamma =0$.
If this condition is not satisfied, we may choose another
imbeding $j': \CClZ \to \X = \Y \oplus \CClZ$, where $j' = j + (\alpha^\dagger)^{-1} \gamma$
(here $j: \CClZ \to \X$ denotes the original imbeding), to achieve $\gamma =0$.

Now we may identify $\X$ with $\Y \oplus \Y$ and assume that the form $\phi$ is given by the
matrix
$$
\bmatrix
0 & \beta\\
\beta & 0
\endbmatrix
$$
where $\beta: \Y \to \e(\Y)$ is a non-degenerate positive form on $\Y$. It is obvious
that we may find a
decomposition $\X = \X_+ \oplus \X_-$ into orthogonal sum of a positive and negative
forms, where $X_+$ is $\Y$, imbeded into $\X$
via $i_1 + i_2$, and $\X_-$ is $\Y$, which is imbeded into $\X$ via $i_1 - i_2$.
The result now follows from Theorem 7.7. \qed

\enddemo


\Refs

\widestnumber\key {GLR}

\ref\key A\by M.F.Atiyah\paper Elliptic operators, discrete groups and von
Neumann algebras\jour Asterisque\vol 32,\-33\yr 1976\pages 43-72\endref

\ref \key BF \by E. Bayer-Fluckiger, L. Fainsilber\paper Non unimodular
Hermitian forms\jour Invent. math. \vol 123\yr 1996\pages 233-240\endref

\ref \key CG\by J.Cheeger and M.Gromov\paper $L_2$-cohomology and
group cohomology\jour Topology\vol 25\yr 1986\pages 189-215
\endref

\ref\key D\by J. A. Dieudonn\'e\paper Quasi-hermitian operators\book
Proc. International Symposium in linear spaces, Jerusalem\yr 1961\pages 115-122
\endref

\ref\key DLS\by A. Dijksma, H. Langer, H. de Snoo\paper Unitary colligations in
Krein spaces and the extension theory of isometries and symmetric linear relations
in Hilbert spaces\jour Lecture Notes in Mathematics\vol 1242\yr 1987\pages 1 - 42
\endref


\ref\key Di\by J.Dixmier\book Von Neumann algebras\yr 1981
\publ North-Holland Publishing Company
\endref

\ref\key DS\by N. Dunford and J. Schwartz\book Linear operators, Part I
\yr 1967 \publ Interscience Publishers\endref

\ref\key{Fa1}\by M. Farber \paper Abelian categories, Novikov-Shubin 
invariants, and Morse inequalities\jour
C. R. Acad. Sci. Paris \vol 321 \yr 1995\pages 1593-1598\endref                       

\ref \key Fa2\by M. Farber\paper Homological algebra of Novikov-Shubin 
invariants and Morse inequalities\jour GAFA\yr 1996\vol 6\endref

\ref\key Fa3\by M. Farber\paper Von Neumann categories and extended
$L^2$ cohomology\jour J. of K-theory (to appear) \yr 1998\endref

\ref\key Fa4\by M. Farber\paper Geometry of growth: approximation theorems for $L^2$
invariants\jour Math. Annalen\yr 1998\vol 311\pages 335 - 375\endref

\ref\key GLR\by P. Ghez, R. Lima, J.E. Roberts\paper $W^\ast$-categories\jour Pacific J. of Math.
\yr 1985\vol 120\pages 79 - 109\endref

\ref\key GK\by I. Gohberg and N. Krupnik\book Introduction to theory of
one-dimensional singular integral operators\yr 1973\endref

\ref\key G\by M.Gromov\book Asymptotic invariants of infinite groups\yr 1993
\publ London Math. Society, Lect. Notes Series\vol 182\endref

\ref\key GS1\by M.Gromov and M.A.Shubin\paper Von Neumann spectra near
zero\jour GAFA\vol 1\yr 1991\pages 375-404\endref

\ref \key GS2\by M.Gromov and M.Shubin\paper Near-cohomology of Hilbert
complexes and topology of non-simp\-ly connected manifolds\jour Asterisque
\vol 210\pages 283-294\yr 1992\endref

\ref\key Gr\by A.Grothendieck
\paper Sur quelques points d'alg\'ebre homologique
\jour Tohoku Mathematical Journal\yr 1957 \vol 9
\endref

\ref\key L\by P. Lax\paper Symmetrizable linear transformations\jour Comm. in Pure and
Applied Math. \vol 7\yr 1954\pages 633 - 647\endref

\ref \key LL\by J.Lott and W.L\"uck\paper $L^2$-topological invariants of
3-manifolds\jour Invent. math. \yr 1995\vol 120\pages 15-60\endref

\ref\key Lu1\by W.L\"uck\paper Hilbert modules and modules over finite von 
Neumann algebras and applications to $L^2$\-in\-va\-riants\jour Math. Annalen
\yr 1997\pages 247 - 285\endref

\ref\key Lu2\by W. L\"uck\paper Dimension theory of arbitrary modules over finite von Neumann 
algebras and $L^2$-Betti numbers I: Foundations\jour J. f\"ur Reine und Angewandte
Mathematik\vol 496\yr 1998\pages 135 - 162\endref

\ref\key Lu3\by W. L\"uck\paper Dimension theory of arbitrary modules over finite von Neumann 
algebras and $L^2$-Betti numbers I: Applications to Grothendieck groups
\jour J. f\"ur Reine und Angewandte
Mathematik\vol 496\yr 1998\pages 213 - 236
\endref

\ref\key Lu4\by W. L\"uck\paper Approximating $L^2$-invariants by their finite dimensional analogues
\jour GAFA\yr 1994\vol 4\pages 455 - 481\endref

\ref\key N\by V.V. Nikulin\paper Integral symmetric bilinear forms and some of their geometrical
applications\jour Math. USSR Izvestija\vol 14\yr 1980\pages 103 - 167\endref

\ref \key NS1\by S.P.Novikov and M.A.Shubin\paper Morse inequalities and
von Neumann invariants of non\-simp\-ly connected manifolds\jour Uspehi
Matem. Nauk \vol 41\yr 1986\pages 222\endref

\ref\key NS2\by S.P.Novikov and M.A.Shubin\paper Morse inequalities and
von Neumann $II_1$-factors\jour Doklady Akad. Nauk SSSR\vol 289\yr 1986
\pages 289-292\endref

\ref \key P\by R.S.Palais\book Seminar on the Atiyah - Singer index theorem
\publ Annals of math. studies, number 57\yr 1965
\endref

\ref  \key QSS\by H.-G. Quebbemann, W. Scharlau, M. Schulte
\paper Quadratic and Hermitian forms in Additive and Abelian Categories
\jour Journal of Algebra\yr 1979\vol 59 \pages 264-289
\endref

\ref \key R1 \by A. Ranicki\paper Algebraic $L$-theory, I: Foundations
\jour Proc. London Math. Soc.\vol 27 \year 1973\pages 101-125
\endref

\ref  \key R2\by A. Ranicki\paper Additive $L$-theory\jour $K$-theory
\vol 3\pages 163-195\yr 1989
\endref


\ref \key W\by C.T.C. Wall \book Surgery on compact manifolds
\publ Academic Press\yr 1970\endref

\ref\key W1\by C.T.C. Wall\paper On the axiomatic foundation of the theory of Hermitian forms
\jour Proc. Cambridge Phil. Soc. \vol 67\yr 1970\pages 243 - 250\endref


\endRefs
\enddocument